\documentclass[11pt]{article}
\usepackage{amsfonts}
\usepackage{amsfonts}
\usepackage{amsfonts}
\usepackage{mathrsfs}
\usepackage{amssymb,amsmath}
 \usepackage{cite}
\begin{document}

\def\pd#1#2{\frac{\partial#1}{\partial#2}}

\let\oldsection\section
\renewcommand\section{\setcounter{equation}{0}\oldsection}
\renewcommand\thesection{\arabic{section}}
\renewcommand\theequation{\thesection.\arabic{equation}}

\newtheorem{thm}{ Theorem}[section]
\newtheorem{lem}{ Lemma}[section]
\newtheorem{proposition}{ Proposition}[section]
\newtheorem{dnt}{ Definition}[section]
\newtheorem{remark}{ Remark}[section]
\newtheorem{cor}{ Corollary}[section]
\allowdisplaybreaks

\title{Strong solutions to the Cauchy problem of the two-dimensional
compressible Navier-Stokes-Smoluchowski equations with vacuum\thanks{This work was
supported by excellent doctorial dissertation cultivation grant from
Dalian
University of
Technology.}}
\author{Yang Liu$\,^{\rm 1}\,${\thanks{Corresponding author.  E-mail address:
liuyang19850524@163.com.}}\\[6pt]
\small $^{\rm 1}\,$School of Mathematical Sciences, Dalian
University of
Technology, Dalian 116024, China \\
\date{}
}\maketitle

\maketitle{}{\bf Abstract.}~   This paper studies the local existence
of strong solutions to the Cauchy problem of the 2D fluid-particle interaction model with
vacuum as far field density. Notice that the technique used by Ding et al.\cite{SBH} for the corresponding
3D local well-posedness of strong solutions fails treating the 2D
case, because the $L^p$-norm ($p>2$) of the velocity $u$ cannot be
controlled in terms only of $\sqrt{\rho}u$ and $\nabla u$
here. In the present paper, we will use the framework of weighted
approximation estimates introduced in [J. Li, Z. Liang, On classical
solutions to the Cauchy problem of the two-dimensional barotropic
compressible Navier-Stokes equations with vacuum, J. Math. Pures
Appl. (2014) 640--671] for Navier-Stokes equations to obtain the
local existence of strong solutions provided the initial density and
 density of
particles in the mixture do not decay very slowly at infinity. In particular, the initial
 density can have a compact support.
 This paper extends Fang et al.'s result \cite{DRZ} and Ding et al.'s result \cite{SBH}, in which, the existence is obtained when the space
dimension $N=1$ and $N=3$ respectively.
\underline{}
\section{Introduction}

In this paper, we consider a fluid-particle interaction model called as
Navier-Stokes-Smoluchowski equations in \cite{JB1, JT1, JT2}, which in the whole spatial domain $\Bbb R^2$ as follows
\begin{equation}\label{a1}
\left\{
\begin{array}{ll}
\displaystyle
 \rho_t+\operatorname{div}(\rho u)=0,\\[3pt]
(\rho u)_t +\operatorname{div}(\rho u\otimes  u)
+\nabla (p_F+\eta)=\mu\Delta
u+(\lambda+\mu)\nabla\operatorname{div}u-(\eta+\beta\rho)\nabla\Phi,\\[3pt]
\eta_t+\nabla\cdot(\eta(u-\nabla\Phi))=\Delta\eta,
\end{array}
\right.
\end{equation}
in $\mathbb{R}^2\times \mathbb{R}^+$, with the far-field behavior
\begin{align}\label{kk3}
(\rho, u, \eta)(x, t)\rightarrow(0, 0, 0)~~ \mathrm{as}~
|x|\rightarrow\infty, ~~ t>0,
\end{align}
and initial data
\begin{equation}\label{a3}
\rho(x, 0)=\rho_0(x), \quad \rho u(x, 0)=m_0, \quad \eta(x, 0)=\eta_0,
\quad x\in \mathbb{R}^2.
\end{equation}
Here $\rho:\mathbb{R}^2\times[0, \infty)\rightarrow\mathbb{R}^+$ is
the density of the fluid, $u:\mathbb{R}^2\times[0,
\infty)\rightarrow \mathbb{R}^2$ the velocity field,
 and the density of the particles in the mixture $\eta:(0, \infty)\times\Bbb R^2\rightarrow\Bbb R_+$
 is related to the probability distribution function $f(t, x, \xi)$ in the
 macroscopic description through the relation
 \begin{align*}
 \eta(t, x)=\int_{\Bbb R^2}f(t, x, \xi)d\xi.
 \end{align*}
 We also denote by $p_F$ the pressure of the fluid, given by
 \begin{align}
 p_F=p_F(\rho)=a\rho^\gamma, a>0, \gamma>1,
 \end{align}
and the time independent external potential $\Phi=\Phi(x):\Bbb R^2\rightarrow\Bbb R_+$
is the effects of gravity and buoyancy, $\beta$ is a constant reflecting the differences in
how the external force affects the fluid and the particles, $\lambda$ and $\mu$ are constant
viscosity coefficients satisfying the physical condition:
\begin{equation}\label{ee}
\mu>0, \quad\lambda+\mu\ge0.\end{equation}

The fluid-particle interaction model arises in a lot of industrial procedures such as the analysis
of sedimentation phenomenon which finds its applications in biotechnology, medicine, chemical
engineering, and mineral processes. Such interaction systems are also used in combustion
theory, when modeling diesel engines or rocket propulsors, see \cite{SR, CL, WS, AF}.
The system  consists in a Vlasov-Fokker-Planck equation to describe the microscopic motion
of the particles coupled to the  equations for the
fluid. Generally speaking, at the microscopic scale, the cloud of particles is described by its distribution function $f(t, x, \xi)$,
solution to a Vlasov-Fokker-Planck equation.
The fluid, on the other hand, is modeled by macroscopic quantities, namely its
density $\rho(x, t)\ge 0$ and its velocity field $u(x, t)$(see \cite{JT1}). If the fluid is compressible and isentropic,
then $(\rho, u)$ solves the compressible Euler (inviscid case) or Navier-Stokes system (viscous case) of
equations.
With the dynamic viscosity terms taken into consideration, system \eqref{a1} was derived formally by Carrillo and Goudon \cite{JT2}.
They obtained the global existence and asymptotic behavior of the weak solutions to \eqref{a1} following
the framework of Lions \cite{PL} and Feireisl et al.\cite{EF1, EF2}. Without the dynamic viscosity
terms in $\eqref{a1}_2$,
Carrillo and Goudon \cite{JT1} gave the
flowing regime
and the bubbling regime under the two different scaling assumptions and investigated the stability and
asymptotic limits finally.
In dimension one, Fang et al.\cite{DRZ}  proved the global existence
and uniqueness of the classical large solution with vacuum.  In dimension three,
 Ballew obtained the local in time existence of strong solutions in a bounded domain
with the no-flux condition for the particle density in \cite{JB1, JB2} and studied Low Mach Number Limits
under the confinement hypotheses for the spatial domain
 and external potential $\Phi$ in \cite{JB3}. Recently, motivated by
 Kim et al. \cite{KIM1, KIM2, KIM3} on the Navier-Stokes equations,
 Ding et al.\cite{SBH} obtained the local classical solutions of
system \eqref{a1} with vacuum in $\Bbb R^3$.

When the density of the
fluid $\eta=0$, the system \eqref{a1} becomes Navier-Stokes equations for the
isentropic compressible
fluids. Kim et al. proved some local existence results on strong solutions in a
domain
 of $\Bbb R^3$ in \cite{KIM1, KIM2} and the radially symmetric solutions in an annular domain in \cite{KIM4}. Ding et
al.\cite{SH} obtained global classical solutions with large initial data with vacuum in a bounded domain or
exterior domain $\Omega$
 of $\Bbb R^n(n\ge 2)$. In a bounded or unbounded domain
 of $\Bbb R^3$, Cho and Kim also got
the local classical solutions \cite{KIM3}, in which the initial density needs not be bounded below away from zero.
For the case that the initial density is allowed
to vanish, Huang et al.\cite{LJ3} obtained the global existence of classical solutions to the
Cauchy problem for the isentropic compressible Navier-Stokes equations in three spatial
dimensions with smooth initial data provided that the initial energy is suitably small.
Recently, assumed that the initial density do not decay very slowly at infinity,
Li and Liang \cite{LJ2} have obtained the local existence of the classical solutions to the two-dimensional
Cauchy problem. After that, Li and Xin \cite{LJ1} extended the result of Li and Liang \cite{LJ2} to the global ones, and also get some
decay estimates of solutions.

The  aim of this paper is to establish the local existence of strong
solutions to the Cauchy problem \eqref{a1} in dimension two.
Notice that the local well-posedness of strong solutions for
dimension three case established by Ding et al.\cite{SBH} is not
admitted for the case of dimension two. This is mainly due to that
in dimension two we fail to control the $L^p$-norm ($p>2$) of the
velocity $u$ in terms only of $\sqrt{\rho}u$ and $\nabla u$.
Moreover, the coupling of $u, \eta$ and $\Phi$, and the presence of $\nabla\cdot(\eta u-\eta\nabla\Phi)$ bring additional difficulties. So, some new ideas and careful
estimates are necessary to deal with the two dimension case. In the
present paper, we will use the framework of weighted approximation
estimates introduced in \cite{LJ2} for Navier-Stokes equations to
overcome these difficulties.
\begin{dnt}
If all derivatives involved in \eqref{a1} for $(\rho, u, \eta)$ are regular distributions,
 and equations \eqref{a1} hold almost everywhere in $\Bbb R^2\times (0, T)$, then $(\rho, u, \eta)$
 is called a strong solution to \eqref{a1}.
\end{dnt}

In this section, for $1\le r\le \infty$, we denote the standard Lebesgue and Sobolev spaces as follows:
\begin{align}
L^r=L^r(\Bbb R^2), \quad W^{s, r}=W^{s, r}(\Bbb R^2), \quad H^s=W^{s, 2}.
\end{align}

 Denote
\begin{equation*}
\bar{x}\triangleq ({\rm e}+|x|^2)^{1/2}\log^{1+\sigma_0}({\rm
e}+|x|^2),
\end{equation*}
with $\sigma_0>0$, $B_N\triangleq\{x\in \mathbb{R}^2 |~ |x|<N\}$. The
main result of this paper is stated as the following theorem:
\begin{thm}\label{thm1}  Suppose that the initial data
$(\rho_0, u_0, \eta_0)$ satisfy
\begin{align*}&\rho_0\ge 0,
\ \bar{x}^a\rho_0\in L^1\cap H^1 \cap W^{1,q}, \ \nabla u_0\in L^2,
 \\ &\nabla \eta_0\in L^2, \
\bar{x}^\frac{a}{2}\eta_0\in L^2, \ \Phi\in H^4, \ \sqrt{\rho}_0u_0\in L^2,
\end{align*}
with $q>2$ and  $a>1$. Then there exist $T_0,N>0$ such that the
problem \eqref{a1}--\eqref{a3} has a unique strong solution
$(\rho, u, \eta)$ on $\mathbb{R}^2\times(0, T_0]$ satisfying
\begin{equation}\label{e7}
\left\{
\begin{array}{ll}
\displaystyle
\rho\in C([0, T_0]; L^1 \cap H^1 \cap W^{1,q}),\quad \bar{x}^a\rho\in L^\infty(0, T_0;  L^1 \cap H^1 \cap W^{1,q}),\\
\sqrt{\rho}u, \nabla u, \bar{x}^{-1}u,\sqrt{t}\sqrt{\rho}u_t\in
L^\infty(0, T_0; L^2),\\
\nabla u\in L^2(0, T_0; H^1)\cap L^{\frac{q+1}{q}}(0, T_0; W^{1,q}),
\sqrt{t}\nabla u\in L^2(0, T_0; W^{1, q}),\\
\eta, \nabla\eta, \bar{x}^{\frac{a}{2}}\eta, \sqrt{t}\eta_t\in L^\infty(0, T_0; L^2),\\
\nabla\eta \in L^2(0, T_0; H^1), \sqrt{t}\nabla u\in L^2(0, T_0,
W^{1, q}),\\
\sqrt{\rho}u_t, \bar{x}^{\frac{a}{2}}\nabla\eta, \sqrt{t}\nabla u_t,
\sqrt{t}\nabla\eta_t, \sqrt{t}\bar{x}^{-1}u_t\in
L^2(\mathbb{R}^2\times(0, T_0)),
\end{array}
\right.
\end{equation}
and
\begin{equation}\label{e8}
\inf_{0\leq t\leq T_0}\int_{B_N}\rho(x,
t)dx\ge\frac{1}{4}\int_{\mathbb{R}^2}\rho_0(x, t)dx.
\end{equation}
\end{thm}

The rest of the paper is organized as follows: In Section 2, we
recall some elementary facts and inequalities used in the sequel.
Sections 3 deals with an approximation problem \eqref{b1} on $B_R$
to derive uniform estimates for the unique strong solution with
respect to $R$. Finally, the proof of Theorem \ref{thm1} will be
given in Section 4.

\section{Preliminaries}

Firstly, the follow local existence theory on bounded ball $B_R\triangleq\{x\in\Bbb R^2: |x|<R\}$, where the initial density
is strictly away from vacuum, can be shown by arguments as in \cite{SBH}.

\begin{lem}\label{mn}~For any given $R>0$ and $B_R=\{x\in \Bbb R^2||x|<R\}$, assume that
$(\rho_0, u_0, \eta_0)$ satisfies
\begin{equation}\label{b2}
\begin{aligned}
&(\rho_0, u_0, \eta_0)\in H^3(B_R), \quad \Phi\in H^4(B_R), \quad \inf_{x\in B_R}\rho_0>0.
\end{aligned}
\end{equation}
Then there exist a small time $T_R>0$ and a unique classical solution $(\rho,
u, \eta)$ to the following initial-boundary-value problem
\begin{equation}\label{b1}
\left\{
\begin{array}{ll}
\displaystyle
 \rho_t+\operatorname{div}(\rho u)=0,\\[3pt]
(\rho u)_t +\operatorname{div}(\rho u\otimes  u)
+\nabla (p_F+\eta)=\mu\Delta
u+(\lambda+\mu)\nabla\operatorname{div}u-(\eta+\beta\rho)\nabla\Phi-R^{-1}u,\\[3pt]
\eta_t+\nabla\cdot(\eta(u-\nabla\Phi))=\Delta\eta,\\[3pt]
u=0, \quad (\nabla\eta+\eta\nabla\Phi)\cdot n=0,\quad x\in\partial B_R, ~t>0,\\[3pt]
(\rho, u, \eta)(x, 0)=(\rho_0, u_0, \eta_0)(x), \quad x\in B_R,
\end{array}
\right.
\end{equation}
on $B_R\times(0, T_R]$ such that
\begin{equation}\label{b3}
\left\{
\begin{array}{ll}
\displaystyle \rho\in C([0, T_R]; H^3),\quad \rho_t\in L^\infty(0,
T_R; H^2),\quad \sqrt{\rho}u_t\in L^\infty(0, T_R; L^2),\\
(u, \eta)\in C([0, T_R]; H^3)\cap L^2(0, T_R; H^4),\quad(u_t, \eta_t)\in L^\infty(0,
T_R; H^1)\cap L^2(0, T_R; H^2),\\
(\sqrt{t}u,\sqrt{t}\eta)\in L^\infty(0, T_R; H^4), \quad(\sqrt{t}u_t, \sqrt{t}\eta_t)
\in L^\infty(0, T_R;H^2),\\
(\sqrt{t}u_{tt}, \sqrt{t}\eta_{tt})\in L^2(0, T_R; H^1), \quad
\sqrt{t}\sqrt{\rho}u_{tt}\in (0, T_R; L^2),\\
 (tu_t, t\eta_t)\in L^\infty(0, T_R; H^3),\quad (tu_{tt}, t\eta_{tt})\in L^\infty(0, T_R; H^1)\cap L^2(0, T_R; H^2),\\
t\sqrt{\rho}u_{ttt}\in L^\infty(0, T_R; L^2),\quad (t^{\frac{3}{2}}u_{tt},
 t^{\frac{3}{2}}\eta_{tt})\in L^\infty(0, T_R; H^2),\\
t^{\frac{3}{2}}\sqrt{\rho}u_{ttt}\in L^\infty(0, T_R; L^2), \quad
(t^{\frac{3}{2}}u_{ttt}, t^{\frac{3}{2}}\eta_{ttt})\in L^2(0, T_R; H^1),\\
\end{array}
\right.
\end{equation}
\end{lem}
where we denote $L^2=L^2(B_R)$ and $H^k=H^k(B_R)$ for positive integer $k$.

Next, for either $\Omega=\Bbb R^2$ or $\Omega=B_R$ with $R\geq1$, the following weighted $L^p$-bounds for
elements of the Hilbert space $\tilde{D}^{1, 2}(\Omega)\triangleq\{v\in H_{{\rm loc}}^1(\Omega)|\nabla v\in L^2(\Omega)\}$
can be found in \cite[theorem B.1]{PL2}.
\begin{lem}
For $m\in[2, \infty)$ and $\theta\in(1+m/2, \infty)$, there exists a positive constant $C$ such that for
either $\Omega=\Bbb R^2$ or $\Omega=B_R$ with $R\geq 1$ and for any $v\in \tilde{D}^{1, 2}(\Omega)$,
\begin{align}
\left(\int_\Omega\frac{|v|^m}{e+|x|^2}(\log(e+|x|^2))^{-\theta}dx\right)^\frac{1}{m}\leq
C\|v\|_{L^2(B_1)}+C\|\nabla v\|_{L^2(\Omega)}.
\end{align}
\end{lem}

\begin{lem}\label{ms}{\rm (Lemma 2.4 in \cite{LJ2})}~ Let $\bar{x}$ and $\sigma_0$ be as in Theorem
\ref{thm1} with $\Omega=\mathbb{R}^2$ or $\Omega=B_{R}$, and
$\rho\in L^1(\Omega)\cap L^\gamma(\Omega)$ with $\gamma>1$  be a
non-negative function satisfying
\begin{equation*}
\int_{B_{N_1}}\rho dx\ge M_1, \quad \int_\Omega \rho^\gamma dx\leq
M_2,
\end{equation*}
with $M_1, M_2>0$, and $B_{N_1}\subset\Omega$ ($N_1\ge1$). Then for
every $v\in \tilde{D}^{1, 2}(\Omega)$, there is $C=C(M_1, M_2, N_1,
\gamma, \sigma_0)>0$ such that
\begin{equation}\label{ks1}
\|v \bar{x}^{-1}\|_{L^2(\Omega)}\leq
C\|\sqrt{\rho}v\|_{L^2(\Omega)}+C\|\nabla v\|_{L^2(\Omega)}.
\end{equation}
 Moreover, for $\varepsilon>0$ and $\sigma>0$ there is $C=C(\varepsilon,
\eta, M_1, M_2, N_1, \gamma, \sigma_0)>0$ such that every $v\in
\tilde{D}^{1, 2}(\Omega)$ satisfies
\begin{equation}\label{ks2}
\|v
\bar{x}^{-\sigma}\|_{L^\frac{2+\varepsilon}{\tilde{\sigma}}(\Omega)}\leq
C\|\sqrt{\rho}v\|_{L^2(\Omega)}+C\|\nabla v\|_{L^2(\Omega)},
\end{equation}
with $\tilde{\sigma}=\min\{1, \sigma\}$.
\end{lem}

Next, the following $L^p$-bound for elliptic systems, whose proof is similar to that of \cite[lemma 12]{KIM1},
is a direct consequence of the combination of a well-known elliptic theory due to
Agmon-Douglis-Nirenberg \cite{SA1, SA2} with a standard scaling procedure.
\begin{lem}
For $p>1$ and $k\geq 0$, there exists a positive constant $C$ depending only on $p$ and $k$
such that
\begin{align}\label{o27}
\|\nabla^{k+2}v\|_{L^p(B_R)}\leq C\|\Delta v\|_{W^{1, p}(B_R)},
\end{align}
for every $v\in W^{k+2, p}(B_R)$ satisfying either
\begin{align*}
v\cdot n=0,~ {\rm rot}v=0, \quad on ~\partial B_R,
\end{align*}
or
\begin{align*}
v=0, \quad on ~\partial B_R.
\end{align*}
\end{lem}

  \section{A priori estimates for approximation problem}
Throughout this section and the next, for $p\in[1, \infty]$ and $k\geq 0$, we denote
\begin{align*}
\int fdx=\int_{B_R}fdx, \quad L^p=L^p(B_R), \quad W^{k, p}=W^{k, p}(B_R), \quad H^k=W^{k, 2}.
\end{align*}
 Moreover, for
$R>4N_0\geq 4$, assume that  $(\rho_0, u_0, \eta_0)$ satisfies, in addition to \eqref{b2}, that
\begin{equation}\label{31}
\frac{1}{2}\leq\int_{B_{N_0}}\rho_0(x)dx\leq\int_{B_R}\rho_0(x)dx\leq \frac{3}{2}.
\end{equation}
Lemma \ref{mn} thus yields that there exists some $T_R>0$ such that the initial-boundary value problem \eqref{b1}
has a unique classical solution $(\rho, u, \eta)$ on $B_R\times[0, T_R]$ satisfying \eqref{b3}.

For $\bar{x}, \sigma_0, a$ and $q$ as in theorem \ref{thm1}, the main aim of this section is to derive the
following key a priori estimate on $\psi$ defined by
\begin{align}
\psi(t)&\triangleq 1+\|\sqrt{\rho}u\|_{L^2}+\|\nabla
u\|_{L^2}+\|\nabla\eta\|_{L^2}
+\|\bar{x}^\frac{a}{2}\eta\|_{L^2}
+\|\bar{x}^a\rho\|_{L^1\cap H^1\cap W^{1, q}}+R^{-1}\|u\|_{L^2}.
\end{align}
\begin{proposition}\label{mk}~ Assume that $(\rho_0, u_0, \eta_0)$
satisfies \eqref{b2} and \eqref{31}. Then there exist $T_0,M>0$,
both depending only on $\mu, \gamma, q, a, \eta_0, N_0$, and $E_0$,
such that
\begin{equation}\label{3.3}
\sup_{0\leq t\leq T_0}\psi(t)+\int_0^{T_0}(\|\nabla^2
u\|_{L^q}^{\frac{q+1}{q}}+t\|\nabla^2 u\|_{L^q}^2+\|\nabla^2
u\|_{L^2}^2+\|\nabla^2\eta\|_{L^2}^2)dt\leq M.
\end{equation}
where
\begin{align*}
E_0&\triangleq\|\sqrt{\rho_0}u_0\|_{L^2}+\|\nabla u_0\|_{L^2}
+\|\nabla\eta_0\|_{L^2}+\|\bar{x}^a\rho_0\|_{L^1\cap H^1\cap W^{1, q}}+\|\bar{x}^\frac{a}{2}\eta_0\|_{L^2},
\end{align*}
\end{proposition}

To prove proposition \ref{mk}, whose proof will be postponed to the end of this section,
we begin with the following standard energy estimate for
$(\rho, u, \eta)$.
\begin{lem}\label{bx1}
Let $(\rho, u, \eta)$ be
a smooth solution to the initial-boundary value problem
\eqref{b1}. Then there exists $T_1=T_1(N_0, E_0)>0$  such that for all $t\in (0, T_1]$
\begin{align}\label{b4}
&\sup_{0\leq s\leq
t}\int\Big[\frac{1}{2}\rho|u|^2+\frac{a}{\gamma-1}\rho^\gamma+\eta\ln\eta
+(\beta\rho+\eta)\Phi\Big]dx\nonumber\\
&\quad+\int_0^t\int\Big[|\nabla u|^2+|2\nabla\sqrt{\eta}+\sqrt{\eta}\nabla\Phi|^2\Big]dxds\leq C,
\end{align}
and moreover,
\begin{align}\label{b5}
\sup_{0\leq s\leq t}\int\eta^2dx+\int_0^t\int|\nabla\eta|^2dxds\leq C.
\end{align}
where and throughout the paper,
denote by $C$ generic positive constants depending only on the fixed constants $\mu, \lambda, \gamma, \beta, a, q, \sigma_0, N_0, E_0$, and $\|\Phi\|_{H^4(\Bbb R^2)}$.

\end{lem}
{\bf Proof.}~
First,
multiplying $\eqref{b1}_2$ by $u$, integrating the resulting equation over $B_R$ and
using Eq.$\eqref{b1}_1$, we have
\begin{align}\label{b6}
&\frac{d}{dt}\int\Big[\frac{1}{2}\rho |u|^2+\frac{a}{\gamma-1}\rho^\gamma\Big]dx\nonumber\\
&+\int\Big[\mu|\nabla u|^2+(\mu+\lambda)|{\rm div}u|^2
+u\cdot\nabla\eta+(\beta\rho+\eta)\nabla\Phi\cdot u+R^{-1}|u|^2\Big]dx=0,
\end{align}
where we have used the fact
\begin{align*}
\int\rho^\gamma\nabla\cdot udx=\int\rho^{\gamma-1}\rho\nabla\cdot udx
=-\int(\rho_t+u\cdot\nabla\rho)\rho^{\gamma-1}=
-\frac{d}{dt}\int\frac{\rho^\gamma}{\gamma}dx+\int\frac{\rho^\gamma}{\gamma}\nabla\cdot u,
\end{align*}
so that
\begin{align*}
-a\int\rho^\gamma\nabla\cdot u=\frac{d}{dt}\int\frac{a\rho^\gamma}{\gamma-1}.
\end{align*}
Using $\eqref{b1}_1$ and $\eqref{b1}_3$, we have
\begin{align}\label{b7}
\int(\eta+\beta\rho)\nabla\Phi\cdot udx
&=-\int{\rm div}(\eta u)dx-\int\beta{\rm div}(\rho u)\Phi dx\nonumber\\
&=\frac{d}{dt}\int\beta\rho\phi dx+\int\Big[\eta_t
-\nabla\cdot(\eta\nabla\phi)-\Delta\eta\Big]\Phi dx\nonumber\\
&=\frac{d}{dt}\int(\eta+\beta\rho)\Phi dx+\int(\eta\nabla\Phi+\nabla\eta)\nabla\Phi dx.
\end{align}
Multiplying $\eqref{b1}_3$ by $\log\eta$, integrating the resulting equation over $B_R$, and using the boundary
condition $\eqref{b1}_4$, one deduces that
\begin{align}\label{b8}
&\int\eta_t\log\eta dx-\int\Big[\eta u-\eta\nabla\Phi-\nabla\eta\Big]\frac{\nabla\eta}{\eta}dx\nonumber\\
&=\frac{d}{dt}\int\eta\log\eta dx-\int\Big[u\cdot\nabla\eta
-\nabla\Phi\cdot\nabla\eta-\frac{|\nabla\eta|^2}{\eta}\Big]dx=0.
\end{align}
Substituting \eqref{b7} and \eqref{b8} into \eqref{b6}, we immediately complete the proof of \eqref{b4}.

Next, multiplying $\eqref{b1}_2$ by $\eta$, integrating the resulting equation over $B_R$, using boundary condition $\eqref{b1}_4$,
we have
\begin{align}
\frac{1}{2}\frac{d}{dt}\int\eta^2dx+\int|\nabla\eta|^2dx&=
\int\eta(u-\nabla\Phi)\nabla\eta dx\nonumber\\
&=\int\eta u\nabla\eta dx-\int\eta\nabla\Phi\nabla\eta dx\nonumber\\
&=-\int{\rm div}u\eta^2dx-\int\eta\nabla\Phi\nabla\eta dx\nonumber\\
&\leq \frac{1}{4}\int|\nabla\eta|^2dx+C\int|\eta|^2|\nabla u|dx+C\int\eta^2|\nabla\Phi|^2dx\nonumber\\
&\leq \frac{1}{4}\int|\nabla\eta|^2dx+C\int\eta^2|\nabla\Phi|^2dx+C\|\nabla u\|_{L^2}\|\eta\|_{L^4}^2\nonumber\\
&\leq \frac{1}{4}\int|\nabla\eta|^2dx+C\int\eta^2|\nabla\Phi|^2dx
+C\|\nabla u\|_{L^2}\|\eta\|_{L^2}\|\nabla\eta\|_{L^2}\nonumber\\
&\leq \frac{1}{2}\int|\nabla\eta|^2dx+C\int\eta^2dx
+C\|\nabla u\|_{L^2}^2\|\eta\|_{L^2}^2.
\end{align}
According to energy inequality \eqref{b4}, we have $\int_0^t\int|\nabla u|^2dxds\leq C$. Thus, we can use Gronwall's inequality to
deduce that
\begin{align}
\sup_{0\leq s\leq t}\int\eta^2dx+\int_0^t\int|\nabla\eta|^2dxds\leq C.
\end{align}

\begin{lem}\label{xx}~
Under the conditions of Proposition \ref{mk}, let $(\rho, u, \eta)$ be
a smooth solution to the initial-boundary value problem
\eqref{b2}-\eqref{b1}. Then there exists $T_1=T_1(N_0, E_0)>0$ and
$\alpha=\alpha(\gamma,q)>1$ such that for all $t\in (0, T_1]$
\begin{align}
&\sup_{0\leq s\leq t}\|\bar{x}^\frac{a}{2}\eta\|_{L^2}^2+\int_0^t\|\bar{x}^\frac{a}{2}\nabla\eta\|_{L^2}^2ds
\leq C,\label{d11}\\
&\sup_{0\leq s\leq t}(\|\nabla u\|_{L^2}^2+\|\nabla\eta\|_{L^2}^2)
+\int_0^t(\|\sqrt{\rho}u_t\|_{L^2}^2+\|\eta_t\|_{L^2}^2+\|\Delta\eta\|_{L^2}^2)ds\leq
 C\int_0^t\psi^\alpha ds+C.\label{c12}
 \end{align}

\end{lem}
{\bf Proof.}~ First, we always assume that
$t\leq T_1$.
The conservation of $\rho$ with $\eqref{b1}_1$ yields that there exists
$T_1>0$ such that
\begin{align}\label{36}
\inf_{0\leq t\leq T_1}\int_{B_{2N_0}}\rho dx\ge \frac{1}{4},
\end{align}
that is (3.8) in \cite{LJ2}. Furthermore, corresponding to (3.10)
obtained in \cite{LJ2},  we have by \eqref{b4}, \eqref{36}, and Lemma \ref{ms} that
\begin{align}\label{3.10}
\|\rho^\eta u\|_{L^{\frac{2+\varepsilon}{\tilde{\sigma}}}}+\|u
\bar{x}^{-\eta}\|_{L^{\frac{2+\varepsilon}{\tilde{\sigma}}}}\leq
C(\varepsilon, \sigma)\psi^{1+\sigma},\quad t\in(0,T_1]
\end{align}
with $\tilde{\sigma}=\min\{1, \sigma\}$.

Next, to obtain \eqref{d11}, multiplying $\eqref{b1}_3$ by
$\eta\bar{x}^a$ and integrating by parts yield
\begin{align}\label{b15}
&\frac{1}{2}\Big(\int|\eta|^2\bar{x}^adx\Big)_t+\int|\nabla\eta|^2\bar{x}^adx\nonumber\\
&=\frac{1}{2}\int\eta^2\Delta\bar{x}^adx-\int\nabla\cdot(\eta u-\eta\nabla\Phi)\eta\bar{x}^adx\nonumber\\
&=\frac{1}{2}\int\eta^2\Delta\bar{x}^adx-\int\eta u\cdot\nabla\eta\bar{x}^adx-\int\eta^2(\nabla\cdot u)\bar{x}^adx\nonumber\\
&\quad+\int\eta\nabla\eta\cdot\nabla\Phi\bar{x}^adx
+\int\eta^2\Delta\Phi\bar{x}^adx\nonumber\\
&\leq C\int|\eta|^2\Delta\bar{x}^adx+C\int|\eta|^2|\nabla u|\bar{x}^adx+C\int|\eta|^2|u|\cdot\nabla\bar{x}^adx\nonumber\\
&\quad+C\int|\eta|^2|\Delta\Phi|\bar{x}^adx
+C\int|\eta|^2|\nabla\Phi|\nabla\bar{x}^adx\nonumber\\
&\triangleq \sum_{i=1}^5I_i.
\end{align}
Direct calculations yield that
\begin{align}
I_1&\leq C\int|\eta|^2\bar{x}^a\bar{x}^{-2}\log^{2(1+\sigma_0)}(e+|x|^2)dx\leq
C\int|\eta|^2\bar{x}^adx,\\
I_2&\leq C\int|\nabla u||\eta|^2\bar{x}^adx\label{b16}\nonumber\\
&\leq C\|\nabla u\|_{L^2}\|\eta\bar{x}^a\|_{L^4}^2\nonumber\\
&\leq C\|\nabla u\|_{L^2}\|\eta\bar{x}^a\|_{L^2}
(\|\nabla\eta\bar{x}^\frac{a}{2}\|_{L^2}+\|\eta\nabla\bar{x}^\frac{a}{2}\|_{L^2})\nonumber\\
&\leq C(\|\nabla u\|_{L^2}^2+1)\|\eta\bar{x}^\frac{a}{2}\|_{L^2}^2+\frac{1}{4}\|\nabla\eta\bar{x}^\frac{a}{2}\|_{L^2}^2,\\
I_3&\leq C\int\bar{x}^a|\eta|^2\bar{x}^{-\frac{3}{4}}|u|\bar{x}^{-\frac{1}{4}}\log^{1+\sigma_0}(e+|x|^2)dx\nonumber\\
&\leq C\|\eta\bar{x}^\frac{a}{2}\|_{L^4}\|\eta\bar{x}^\frac{a}{2}\|_{L^2}\|u\bar{x}^{-\frac{3}{4}}\|_{L^4}\nonumber\\
&\leq C\|\eta\bar{x}^\frac{a}{2}\|_{L^4}^2+C\|\eta\bar{x}^\frac{a}{2}\|_{L^2}^2(\|\sqrt{\rho}u\|_{L^2}^2
+\|\nabla u\|_{L^2}^2)\nonumber\\
&\leq C(1+\|\nabla u\|_{L^2}^2)\|\eta\bar{x}^\frac{a}{2}\|_{L^2}^2+\frac{1}{4}\|\nabla\eta\bar{x}^\frac{a}{2}\|_{L^2}^2,\\
I_4+I_5&\leq C\int|\eta|^2\bar{x}^adx+C\int|\eta|^2\bar{x}^a\bar{x}^{-1}\log^{1+\sigma_0}(e+|x|^2)dx\nonumber\\
&\leq C\int|\eta|^2\bar{x}^adx.\label{b19}
\end{align}
Putting \eqref{b16}-\eqref{b19} into \eqref{b15}, after using Gronwall's inequality and \eqref{b4}, we have
\begin{align}\label{o20}
\sup_{0\leq s\leq t}\|\bar{x}^\frac{a}{2}\eta\|_{L^2}^2+\int_0^t\|\bar{x}^\frac{a}{2}\nabla\eta\|_{L^2}^2dx
\leq C\exp\Big\{C\int_0^t(1+\|\nabla u\|_{L^2}^2)ds\Big\}\leq C.
\end{align}

Next, to prove \eqref{c12},  multiplying Eqs. $\eqref{b1}_2$ by $u_t$ and integration by parts yield
\begin{align}\label{p27}
&\frac{1}{2}\frac{d}{dt}\int\Big[(2\mu+\lambda)({\rm div}u)^2
+\mu\omega^2+R^{-1}|u|^2\Big]dx+\int\rho|u_t|^2dx\nonumber\\
&\leq C\int\rho|u|^2|\nabla u|^2dx
+2\int(p_F+\eta){\rm div}u_tdx-\int(\beta\rho+\eta)\nabla\Phi\cdot u_tdx,
\end{align}
where $\omega\triangleq{\rm rot} u$ is defined in the following \eqref{p44}.

We estimate each term on the right-hand side of \eqref{p27} as follows:

First, the Gagliardo-Nirenberg inequality implies that for all $p\in (2, +\infty)$,
\begin{align}\label{p28}
\|\nabla u\|_{L^p}&\leq C(p)\|\nabla u\|_{L^2}^{2/p}\|\nabla u\|_{H^1}^{1-2/p}\nonumber\\
&\leq C(p)\psi+C(p)\psi\|\nabla^2 u\|_{L^2}^{1-2/p},
\end{align}
which together with \eqref{3.10} yields that for $\sigma>0$ and $\tilde{\sigma}=\min\{1, \sigma\}$,
\begin{align}\label{p29}
\int\rho^\sigma|u|^2|\nabla u|^2dx&\leq
C\|\rho^{\sigma/2}u\|_{L^{8/\tilde{\sigma}}}^2\|\nabla u\|_{L^{8/(4-\tilde{\sigma})}}^2\nonumber\\
&\leq C(\sigma)\psi^{4+2\sigma}(1+\|\nabla^2u\|_{L^2}^{\tilde{\sigma}/2})\nonumber\\
&\leq C(\varepsilon, \sigma)\psi^{\alpha(\sigma)}+\varepsilon\psi^{-2}\|\nabla^2u\|_{L^2}^2.
\end{align}

Next, noticing that $p_F$ satisfies
\begin{align}\label{p30}
p_{Ft}+{\rm div}(p_F u)+(\gamma-1)p_F{\rm div}u=0.
\end{align}
we deduce from $\eqref{b1}_1$ and the Sobolev inequality that
\begin{align}
2\int p_F{\rm div}u_tdx&=2\frac{d}{dt}\int p_F{\rm div}udx-2\int p_F u\cdot\nabla{\rm div}udx
+2(\gamma-1)\int p_F({\rm div}u)^2dx\nonumber\\
&\leq 2\frac{d}{dt}\int p_F{\rm div}udx+\varepsilon\psi^{-1}\|\nabla^2 u\|_{L^2}^2+C(\varepsilon)\psi^\alpha.
\end{align}

Moreover, we have
\begin{align}\label{p32}
&-\int(\nabla\eta+\eta\nabla\Phi)\cdot u_tdx-\int\beta\rho\nabla\Phi\cdot u_tdx\nonumber\\
&=\frac{d}{dt}\int(\nabla\eta+\eta\nabla\Phi)\cdot udx
-\int(\nabla\eta+\eta\nabla\Phi)_t\cdot udx-\int\beta\rho\nabla\Phi\cdot u_tdx\nonumber\\
&=\frac{d}{dt}\int(\eta\nabla\cdot u-\eta\nabla\Phi \cdot u)dx-\int\eta_t\nabla\cdot udx
+\int\eta_t\nabla\Phi \cdot udx-\int\beta\rho\nabla\Phi\cdot u_tdx\nonumber\\
&=\frac{d}{dt}\int(\eta\nabla\cdot u-\eta\nabla\Phi \cdot u)dx-\int\eta_t\nabla\cdot udx\nonumber\\
&\quad+\int(\Delta\eta-\nabla\cdot(\eta u-\eta\nabla\Phi))\nabla\Phi \cdot udx-\int\beta\rho\nabla\Phi\cdot u_tdx\nonumber\\
&= \frac{d}{dt}\int(\eta\nabla\cdot u-\eta\nabla\Phi \cdot u)dx-\int\eta_t\nabla\cdot udx\nonumber\\
&\quad-\int\nabla u\cdot\nabla\eta\cdot\nabla\Phi dx-\int u\cdot\nabla\eta\cdot\Delta\Phi dx
+\int\eta u\cdot\nabla u\cdot\nabla\Phi dx-\int\beta\rho\nabla\Phi\cdot u_tdx\nonumber\\
&\quad+\int\eta u^2\cdot\Delta\Phi dx
-\int\eta u\cdot\nabla\Phi\cdot\Delta\Phi dx-\int\eta\nabla u|\nabla\Phi|^2 dx,\nonumber\\
&\triangleq \frac{d}{dt}J_0+\sum_{i=1}^8J_i.
\end{align}
Direct calculations yield that
\begin{align}
J_1&\leq C\int|\eta_t||\nabla u|dx\leq \frac{1}{2}\|\eta_t\|_{L^2}^2+C\int|\nabla u|^2dx,\label{p33}\\
J_2&\leq C\|\nabla\Phi\|_{L^\infty}\int|\nabla u||\nabla\eta|dx\leq C\|\nabla\eta\|_{L^2}^2+C\|\nabla u\|_{L^2}^2,\\
J_3&\leq C\int|u|\bar{x}^{-\frac{a}{2}}|\nabla\eta|\bar{x}^\frac{a}{2}\Delta\Phi dx\nonumber\\
&\leq C\|\Delta\Phi\|_{L^4}\|u\bar{x}^{-\frac{a}{2}}\|_{L^4}\|\nabla\eta\bar{x}^\frac{a}{2}\|_{L^2}\nonumber\\
&\leq C\|\bar{x}^\frac{a}{2}\nabla\eta\|_{L^2}^2+C(1+\|\nabla u\|_{L^2}^2),\\
J_4&\leq C\int|\eta|\bar{x}^\frac{a}{2}|u|\bar{x}^{-\frac{a}{2}}|\nabla u||\nabla\Phi|dx\nonumber\\
&\leq C\|\nabla\Phi\|_{L^\infty}\|\eta\bar{x}^\frac{a}{2}\|_{L^2}\|u\bar{x}^{-\frac{a}{2}}\|_{L^4}\|\nabla u\|_{L^4}\nonumber\\
&\leq C\|\eta\bar{x}^\frac{a}{2}\|_{L^2}^2+C(1+\|\nabla u\|_{L^2}^2)\|\nabla u\|_{L^2}\|\nabla^2 u\|_{L^2}\nonumber\\
&\leq C\psi^\alpha+\varepsilon\psi^{-1}\|\nabla^2 u\|_{L^2}^2,\\
J_5&\leq C\|\nabla\Phi\|_{L^\infty}\int\sqrt{\rho}\sqrt{\rho} u_tdx\nonumber\\
&\leq \frac{1}{2}\int\rho|u_t|^2dx+C\|\nabla\Phi\|_{L^\infty}\int\rho_0dx\nonumber\\
&\leq \frac{1}{2}\int\rho|u_t|^2dx+C,\\
J_6&\leq C\int|\eta|\bar{x}^\frac{a}{2}|u|^2\bar{x}^{-\frac{a}{2}}|\Delta\Phi|dx\nonumber\\
&\leq C\|\Delta\Phi\|_{L^\infty}\|\eta\bar{x}^\frac{a}{2}\|_{L^2}\|u\bar{x}^{-\frac{a}{4}}\|_{L^4}^2\nonumber\\
&\leq C\|\Delta\Phi\|_{L^\infty}\|\eta\bar{x}^\frac{a}{2}\|_{L^2}
(\|\sqrt{\rho} u\|_{L^2}^2+\|\nabla u\|_{L^2}^2)\nonumber\\
&\leq C\psi^\alpha,\\
J_7&\leq C\|\Delta\Phi\|_{L^\infty}\int|\eta|\bar{x}^\frac{a}{2}|u|\bar{x}^{-\frac{a}{2}}|\nabla\Phi|dx\nonumber\\
&\leq C\|\eta\bar{x}^\frac{a}{2}\|_{L^4}\|u\bar{x}^{-\frac{a}{2}}\|_{L^4}\|\nabla\Phi\|_{L^2}\nonumber\\
&\leq C\|\bar{x}^\frac{a}{2}\nabla\eta\|_{L^2}^2+C\psi^\alpha,\\
J_8&\leq C\|\nabla\Phi\|_{L^\infty}^2\int|\eta||\nabla u|dx\nonumber\\
&\leq C\|\nabla\Phi\|_{L^\infty}^2\|\eta\|_{L^2}\|\nabla u\|_{L^2}\nonumber\\
&\leq C\psi^\alpha+C.\label{p44}
\end{align}
Substituting \eqref{p33}-\eqref{p44} into \eqref{p32}, and combining \eqref{p32} and \eqref{p27} lead to
\begin{align}\label{p41}
&\frac{1}{2}\frac{d}{dt}\int\Big[(2\mu+\lambda)({\rm div}u)^2
+\mu\omega^2+R^{-1}|u|^2\Big]dx+\int\rho|u_t|^2dx\nonumber\\
&\leq \frac{d}{dt}B(t)+\varepsilon\psi^{-1}\|\nabla^2 u\|_{L^2}^2
+C\|\bar{x}^\frac{a}{2}\nabla\eta\|_{L^2}^2+C\psi^\alpha,
\end{align}
where
\begin{align}
B(t)&=-2\int p_F{\rm div}udx+\int(\eta\nabla\cdot u-\eta\nabla\Phi u)dx\nonumber\\
&\leq \frac{\mu}{4}\|\nabla u\|_{L^2}^2+C\|p_F\|_{L^2}^2+C\|\eta\|_{L^2}^2
+C\|\nabla\Phi\|_{L^4}\|\bar{x}^\frac{a}{2}\eta\|_{L^2}\|u\bar{x}^{-\frac{a}{2}}\|_{L^4}\nonumber\\
&\leq \frac{\mu}{2}\|\nabla u\|_{L^2}^2+C\|p_F\|_{L^2}^2+C,
\end{align}
owing to \eqref{b4}, \eqref{b5}, \eqref{o20} and \eqref{3.10}.

 Moreover, multiplying the equation $\eqref{b1}_3$ by $\eta_t$ and integrating the result equation with respect to $x$ over
$B_R$, we have
\begin{align}\label{b13}
&\frac{1}{2}\frac{d}{dt}\|\nabla\eta\|_{L^2}^2+\|\eta_t\|_{L^2}^2\nonumber\\
&=\int\nabla\cdot(\eta(u-\nabla\Phi))\cdot\eta_tdx\nonumber\\
&\leq\frac{1}{2}\|\eta_t\|_{L^2}^2+C\int |u|^2|\nabla\eta|^2dx
+C\int|\eta|^2|{\rm div}u|^2dx+C\int|\nabla\eta|^2|\nabla\Phi|^2dx+C\int|\eta|^2|\Delta\Phi|^2dx\nonumber\\
&\leq \frac{1}{2}\|\eta_t\|_{L^2}^2+C\int|u|^2\bar{x}^{-\frac{a}{2}}|\nabla\eta|\bar{x}^\frac{a}{2}|\nabla\eta|dx
+C\|\eta\|_{L^4}^2\|\nabla u\|_{L^4}^2+C\int|\nabla\eta|^2dx+C\nonumber\\
&\leq \frac{1}{2}\|\eta_t\|_{L^2}^2+C\|\bar{x}^{-\frac{a}{4}}u\|_{L^8}^2\|\bar{x}^\frac{a}{2}\nabla\eta\|_{L^2}\|\nabla\eta\|_{L^4}
+C\|\eta\|_{L^2}\|\nabla\eta\|_{L^2}\|\nabla u\|_{L^2}\|\nabla u\|_{H^1}+C\int|\nabla\eta|^2dx+C\nonumber\\
&\leq \frac{1}{2}\|\eta_t\|_{L^2}^2+\varepsilon\psi^{-1}\|\nabla^2 u\|_{L^2}^2
+C\|\bar{x}^\frac{a}{2}\nabla\eta\|_{L^2}^2+C\psi^\alpha\|\nabla\eta\|_{L^4}^2
+C(\varepsilon)\psi^\alpha\nonumber\\
&\leq \frac{1}{2}\|\eta_t\|_{L^2}^2+\frac{1}{2}\|\nabla^2\eta\|_{L^2}^2+\varepsilon\psi^{-1}\|\nabla^2 u\|_{L^2}^2
+C\|\bar{x}^\frac{a}{2}\nabla\eta\|_{L^2}^2+C(\varepsilon)\psi^\alpha.
\end{align}
From $\eqref{b1}_3$, taking it by $L^2$-norm, using Gagliardo-Nirenberg inequality, we get
\begin{align}\label{b14}
\|\Delta\eta\|_{L^2}&\leq C\|\eta_t\|_{L^2}+C\|\nabla\cdot(\eta u-\eta\nabla\Phi)\|_{L^2}\nonumber\\
&\leq C\|\eta_t\|_{L^2}+C(\|u\nabla\eta\|_{L^2}+\|\eta\nabla u\|_{L^2}
+\|\nabla\Phi\|_{L^\infty}\|\nabla\eta\|_{L^2}+\|\Delta\Phi\|_{L^\infty}\|\eta\|_{L^2})\nonumber\\
&\leq C\|\eta_t\|_{L^2}+C\|\bar{x}^{-\frac{a}{4}}u\|_{L^8}\|\bar{x}^\frac{a}{2}\nabla\eta\|_{L^2}^\frac{1}{2}\|\nabla\eta\|_{L^4}^\frac{1}{2}\nonumber\\
&\quad+C\|\eta\|_{L^2}^\frac{1}{2}\|\nabla\eta\|_{L^2}^\frac{1}{2}\|\nabla u\|_{L^2}^\frac{1}{2}\|\nabla u\|_{H^1}^\frac{1}{2}
+C\|\nabla\eta\|_{L^2}+C\nonumber\\
&\leq C\|\eta_t\|_{L^2}+\frac{1}{2}\|\nabla^2\eta\|_{L^2}+\varepsilon\psi^{-1}\|\nabla^2 u\|_{L^2}
+C\|\bar{x}^\frac{a}{2}\nabla\eta\|_{L^2}+C(\varepsilon)\psi^\alpha
\end{align}

Finally, to estimate the last term on the right-hand side of \eqref{b13}, \eqref{b14}, and \eqref{p41},
denoting $\nabla^\bot\triangleq(\partial_2, -\partial_1)$,
we rewrite the momentum equation $\eqref{b1}_2$ as
\begin{align}\label{p43}
R^{-1}u+\rho\dot{u}=\nabla F+\mu\nabla^\bot\omega-(\eta+\beta\rho)\nabla\Phi,
\end{align}
where
\begin{align}\label{p44}
\dot{f}\triangleq f_t+u\cdot\nabla f, \quad F\triangleq(2\mu+\lambda){\rm div}u-p_F-\eta,
\quad \omega\triangleq\nabla^\bot\cdot u
\end{align}
are the material derivative of $f$, the effective viscous flux and the vorticity respectively.
Thus, \eqref{p43} implies that $\omega$ satisfies
\begin{equation}\label{p45}
\left\{
\begin{array}{ll}
\displaystyle
\mu\Delta\omega=\nabla^\bot\big(\rho\dot{u}+(\eta+\beta\rho)\nabla\Phi+R^{-1}u\big), &{\rm in}~B_R,\\
\omega=0, &{\rm on}~\partial B_R.
\end{array}
\right.
\end{equation}
Applying the standard $L^p$-estimate to \eqref{p45} yields that, for $p\in (1, \infty)$,
\begin{align*}
\|\nabla\omega\|_{L^p}\leq C(\rho)(\|\rho\dot{u}\|_{L^p}
+\|(\eta+\beta\rho)\nabla\Phi\|_{L^p}+R^{-1}\|u\|_{L^p}),
\end{align*}
which together with \eqref{p43} gives
\begin{align}\label{p46}
\|\nabla F\|_{L^p}+\|\nabla\omega\|_{L^p}\leq C(\rho)(\|\rho\dot{u}\|_{L^p}
+\|(\eta+\beta\rho)\nabla\Phi\|_{L^p}+R^{-1}\|u\|_{L^p}),
\end{align}
It follows from \eqref{o27} and \eqref{p46} that for $p\in [2, q]$,
\begin{align}\label{p47}
\|\nabla^2 u\|_{L^p}&\leq C\|\nabla\omega\|_{L^p}+C\|\nabla{\rm div}u\|_{L^p}\nonumber\\
&\leq C(\|\nabla\omega\|_{L^p}+\|(2\mu+\lambda)\nabla{\rm div}u\|_{L^p}
)\nonumber\\
&\leq C(\|\rho\dot{u}\|_{L^p}
+\|(\eta+\beta\rho)\nabla\Phi\|_{L^p}
+\|\nabla p_F\|_{L^p}+\|\nabla\eta\|_{L^p}+R^{-1}\|u\|_{L^p}),
\end{align}
which together with \eqref{b5}, \eqref{p28} and \eqref{p29} leads to
\begin{align}\label{p48}
\|\nabla^2 u\|_{L^2}&\leq C\sqrt{\psi}\|\sqrt{\rho}u_t\|_{L^2}+C\|\rho u\cdot\nabla u\|_{L^2}
+C\psi^\alpha\nonumber\\
&\leq C\sqrt{\psi}\|\sqrt{\rho}u_t\|_{L^2}+C\psi^\alpha+\frac{1}{2}\|\nabla^2 u\|_{L^2}.
\end{align}
Putting \eqref{p48} into \eqref{b13}, \eqref{b14}, and \eqref{p41}, integrating the resulting inequality over $(0, t)$ and
choosing $\varepsilon$ suitably small yield
\begin{align}
&R^{-1}\|u\|_{L^2}^2+\|\nabla u\|_{L^2}^2+\|\nabla\eta\|_{L^2}^2
+\int_0^t(\|\sqrt{\rho}u_t\|_{L^2}^2+\|\eta_t\|_{L^2}^2+\|\nabla^2\eta\|_{L^2}^2)ds\nonumber\\
&\leq C+C\|p_F\|_{L^2}^2+C\int_0^t\psi^\alpha ds
+C\int_0^t\|\bar{x}^\frac{a}{2}\nabla\eta\|_{L^2}^2ds\nonumber\\
&\leq C+C\|p_F\|_{L^2}^2+C\int_0^t\psi^\alpha ds\nonumber\\
&\leq C+C\int_0^t\psi^\alpha ds,
\end{align}
where we have used \eqref{b15} and the following estimate
\begin{align}
\|p_F\|_{L^2}^2\leq \|p_F(\rho_0)\|_{L^2}^2
+C\int_0^t\|p_F\|_{L^1}^{1/2}\|p_F\|_{L^\infty}^{3/2}\|\nabla u\|_{L^2}ds\leq C+C\int_0^t\psi^\alpha ds
\end{align}
due to \eqref{p30}. The proof of lemma \ref{xx} is completed.
\begin{lem}\label{lm3}
Let $(\rho, u, \eta)$ and $T_1$ be as in Lemma \ref{xx}. Then, for all $t\in (0, T_1]$,
\begin{align}
&\sup_{0\leq s\leq t}s\|\bar{x}^\frac{a}{2}\nabla\eta\|_{L^2}^2
+\int_0^ts\|\bar{x}^\frac{a}{2}\Delta\eta\|_{L^2}^2ds\leq C\exp\Big\{\int_0^t\psi^\alpha ds\Big\},\\
&\sup_{0\leq s\leq t}s\|\sqrt{\rho}u_t\|_{L^2}^2
+\int_0^ts\big(\|\nabla u_t\|_{L^2}^2+R^{-1}|u_t|^2\big)ds
\leq C\exp\Big\{C\int_0^t\psi^\alpha ds\Big\}\label{p52}.
\end{align}
\end{lem}
{\bf Proof. }~
Differentiating $\eqref{b1}_2$ with respect to $t$ gives
\begin{align}\label{c53}
&\rho u_{tt}+\rho u\cdot\nabla u_t-\mu\nabla^\bot\omega_t-\nabla((2\mu+\lambda){\rm div}u_t)+R^{-1}u_t\nonumber\\
&=-\rho_t(u_t+u\cdot\nabla u)-\rho u_t\cdot\nabla u-\nabla(p_{Ft}+\eta_t)
-(\eta_t+\beta\rho_t)\nabla\Phi.
\end{align}
Multiplying \eqref{c53} by $u_t$ and integrating the resulting equation over $B_R$, we obtain after using $\eqref{b1}_1$ that
\begin{align}\label{c54}
&\frac{1}{2}\frac{d}{dt}\int\rho|u_t|^2dx
+\int\big((2\mu+\lambda)({\rm div}u_t)^2+\mu\omega_t^2+R^{-1}|u_t|^2\big)dx\nonumber\\
&=-2\int\rho u\cdot\nabla u_t\cdot u_tdx-\int\rho u\cdot\nabla(u\cdot\nabla u\cdot u_t)dx
\nonumber\\
&\quad-\int\rho u_t\cdot\nabla u\cdot u_tdx+\int (p_{Ft}+\eta_t){\rm div}u_tdx
-\int(\eta_t+\beta\rho_t)\nabla\Phi \cdot u_tdx\nonumber\\
&\triangleq \Psi(t)+\int\eta_t{\rm div}u_tdx-\int(\eta_t+\beta\rho_t)\nabla\Phi \cdot u_tdx.
\end{align}
 By the arguments (3.27)--(3.31) for the proof of Lemma 3.3 in \cite{LJ2}, it follows from \eqref{b4}, \eqref{b5}
and \eqref{3.10} for $\varepsilon\in (0, 1)$ that
\begin{align}
\Psi(t)\leq \varepsilon\|\nabla u_t\|_{L^2}^2+C(\varepsilon)\psi^\alpha(\|\nabla^2 u\|_{L^2}^2
+\|\rho^{1/2}u_t\|_{L^2}^2+1).
\end{align}
On the other hand,
\begin{align}
&\int\eta_t{\rm div}u_tdx-\int(\eta_t+\beta\rho_t)\nabla\Phi u_tdx\nonumber\\
&= \int(\Delta\eta-\nabla\cdot(\eta(u-\nabla\Phi)))\cdot{\rm div}u_tdx
-\int(\Delta\eta-\nabla\cdot(\eta(u-\nabla\Phi))-\beta{\rm div}(\rho u))\nabla\Phi\cdot u_tdx\nonumber\\
&\leq C\int|\Delta\eta||\nabla u_t|dx+C\int|u||\nabla\eta||\nabla u_t|dx
+C\int\eta|\nabla u||\nabla u_t|dx\nonumber\\
&\quad+C\int|\nabla\eta||\nabla u_t||\nabla\Phi|dx+C\int\eta|\nabla u_t||\Delta\Phi|dx
+C\int|\nabla\eta||\Delta\Phi||u_t|dx\nonumber\\
&\quad+C\int\eta|u||\Delta\Phi||u_t|dx+C\int\eta|u||\nabla\Phi||\nabla u_t|dx
+C\int\eta|\nabla\Phi||\Delta\Phi||u_t|dx\nonumber\\
&\quad+C\int\eta|\nabla\Phi|^2|\nabla u_t|dx
+C\int|\nabla\rho||u||\nabla\Phi||u_t|dx+C\int\rho|\nabla u||\nabla\Phi||u_t|dx\nonumber\\
&\triangleq\sum_{i=1}^{12}R_i.
\end{align}
Using Gagliardo-Nirenberg and H\"older's inequalities, we get
\begin{align}
R_1&\leq C\|\Delta\eta\|_{L^2}\|\nabla u_t\|_{L^2}\nonumber\\
&\leq \frac{1}{12}\|\nabla u_t\|_{L^2}^2+C\|\eta_t\|_{L^2}^2+\varepsilon\psi^{-1}\|\nabla^2 u\|_{L^2}^2
+C\|\bar{x}^\frac{a}{2}\nabla\eta\|_{L^2}^2+C(\varepsilon)\psi^\alpha,\label{c57}\\
R_2&\leq C\int|u|\bar{x}^{-\frac{a}{2}}|\nabla\eta|\bar{x}^\frac{a}{2}|\nabla u_t|dx\nonumber\\
&\leq C\|u\bar{x}^{-\frac{a}{2}}\|_{L^4}\|\bar{x}^\frac{a}{2}\nabla\eta\|_{L^4}\|\nabla u_t\|_{L^2}\nonumber\\
&\leq \frac{1}{12}\|\nabla u_t\|_{L^2}^2+C\|u\bar{x}^{-\frac{a}{2}}\|_{L^4}^2
\|\bar{x}^\frac{a}{2}\nabla\eta\|_{L^2}\|\bar{x}^\frac{a}{2}\nabla^2\eta\|_{L^2}\nonumber\\
&\leq\frac{1}{12}\|\nabla u_t\|_{L^2}^2+\frac{1}{8}\|\bar{x}^\frac{a}{2}\nabla^2\eta\|_{L^2}^2
+C\psi^\alpha\|\bar{x}^\frac{a}{2}\nabla\eta\|_{L^2}^2,\\
R_3&\leq C\|\nabla u_t\|_{L^2}\|\eta\|_{L^4}\|\nabla u\|_{L^4}\nonumber\\
&\leq C\|\nabla u_t\|_{L^2}\|\nabla\eta\|_{L^2}^\frac{1}{2}\|\eta\|_{L^2}^\frac{1}{2}
\|\nabla u\|_{L^2}^\frac{1}{2}\|\nabla^2u\|_{L^2}^\frac{1}{2}\nonumber\\
&\leq \frac{1}{12}\|\nabla u_t\|_{L^2}^2+\varepsilon\psi^{-1}\|\nabla^2 u\|_{L^2}^2
+C(\varepsilon)\psi^\alpha,\\
R_4&\leq C\|\nabla\Phi\|_{L^\infty}\|\nabla\eta\|_{L^2}\|\nabla u_t\|_{L^2}\leq \frac{1}{12}\|\nabla u_t\|_{L^2}^2+C(\varepsilon)\psi^\alpha,\\
R_5&\leq C\|\Delta\Phi\|_{L^\infty}\|\eta\|_{L^2}\|\nabla u_t\|_{L^2}\leq \frac{1}{12}\|\nabla u_t\|_{L^2}^2+C,\\
R_6&\leq C\int\bar{x}^\frac{a}{2}|\nabla\eta|\bar{x}^{-\frac{a}{2}}|u_t||\Delta\Phi|dx\nonumber\\
&\leq C\|\bar{x}^\frac{a}{2}\nabla\eta\|_{L^4}\|u_t\bar{x}^{-\frac{a}{2}}\|_{L^4}\|\Delta\Phi\|_{L^2}\nonumber\\
&\leq C\|\bar{x}^\frac{a}{2}\nabla\eta\|_{L^2}^\frac{1}{2}\|\bar{x}^\frac{a}{2}\Delta\eta\|_{L^2}^\frac{1}{2}
(\|\sqrt{\rho}u_t\|_{L^2}+\|\nabla u_t\|_{L^2})\nonumber\\
&\leq  \frac{1}{12}\|\nabla u_t\|_{L^2}^2+C\|\bar{x}^\frac{a}{2}\nabla\eta\|_{L^2}\|\bar{x}^\frac{a}{2}\Delta\eta\|_{L^2}
+C\|\sqrt{\rho}u_t\|_{L^2}^2,\\
R_7&\leq C\int\bar{x}^\frac{a}{2}\eta\bar{x}^{-\frac{a}{4}}|u||\Delta\Phi|\bar{x}^{-\frac{a}{4}}|u_t|dx\nonumber\\
&\leq C\|\Delta\Phi\|_{L^\infty}\|\bar{x}^\frac{a}{2}\eta\|_{L^2}\|u\bar{x}^{-\frac{a}{4}}\|_{L^4}
\|u_t\bar{x}^{-\frac{a}{4}}\|_{L^4}\nonumber\\
&\leq C\|\bar{x}^\frac{a}{2}\eta\|_{L^2}(\|\sqrt{\rho}u\|_{L^2}+\|\nabla u\|_{L^2})
(\|\sqrt{\rho}u_t\|_{L^2}+\|\nabla u_t\|_{L^2})\nonumber\\
&\leq \frac{1}{12}\|\nabla u_t\|_{L^2}^2+C\psi^\alpha\|\sqrt{\rho}u_t\|_{L^2}^2+C\psi^\alpha,\\
R_8&\leq  C\int\bar{x}^\frac{a}{2}\eta\bar{x}^{-\frac{a}{2}}|u||\nabla\Phi||\nabla u_t|dx\nonumber\\
&\leq C\|\nabla\Phi\|_{L^\infty}\|\bar{x}^\frac{a}{2}\eta\|_{L^4}\|u\bar{x}^{-\frac{a}{2}}\|_{L^4}\|\nabla u_t\|_{L^2}\nonumber\\
&\leq C\|\bar{x}^\frac{a}{2}\eta\|_{L^2}^\frac{1}{2}\|\bar{x}^\frac{a}{2}\nabla\eta\|_{L^2}^{1/2}
(\|\sqrt{\rho}u\|_{L^2}+\|\nabla u\|_{L^2})\|\nabla u_t\|_{L^2}\nonumber\\
&\leq \frac{1}{12}\|\nabla u_t\|_{L^2}^2+C\|\bar{x}^\frac{a}{2}\nabla\eta\|_{L^2}^2+C\psi^\alpha,\\
R_9&\leq C\int\bar{x}^\frac{a}{2}\eta|\nabla\Phi||\Delta\Phi|\bar{x}^{-\frac{a}{2}}|u_t|dx\nonumber\\
&\leq C\|\Delta\Phi\|_{L^\infty}\|\nabla\Phi\|_{L^2}\|\bar{x}^\frac{a}{2}\eta\|_{L^4}\|u_t\bar{x}^{-\frac{a}{2}}\|_{L^4}\nonumber\\
&\leq C\|\bar{x}^\frac{a}{2}\eta\|_{L^2}^\frac{1}{2}\|\bar{x}^\frac{a}{2}\nabla\eta\|_{L^2}^\frac{1}{2}
(\|\sqrt{\rho}u_t\|_{L^2}+\|\nabla u_t\|_{L^2})\nonumber\\
&\leq \frac{1}{12}\|\nabla u_t\|_{L^2}^2+C\|\bar{x}^\frac{a}{2}\nabla\eta\|_{L^2}^2
+C\psi^\alpha\|\sqrt{\rho}u_t\|_{L^2}^2+C\psi^\alpha,\nonumber\\
R_{10}&\leq C\|\nabla\Phi\|_{L^\infty}^2\|\eta\|_{L^2}\|\nabla u_t\|_{L^2}\leq \frac{1}{12}\|\nabla u_t\|_{L^2}^2+C,\\
R_{11}&\leq C\int\bar{x}^\frac{a}{2}|\nabla\rho|\bar{x}^{-\frac{a}{4}}|u||\nabla\Phi|\bar{x}^{-\frac{a}{4}}|u_t|dx\nonumber\\
&\leq C\|\nabla\Phi\|_{L^\infty}\|\bar{x}^\frac{a}{2}\nabla\rho\|_{L^2}\|u\bar{x}^{-\frac{a}{4}}\|_{L^4}
\|u_t\bar{x}^{-\frac{a}{4}}\|_{L^4}\nonumber\\
&\leq C\psi^\alpha(\|\sqrt{\rho}u\|_{L^2}+\|\nabla u\|_{L^2})(\|\sqrt{\rho}u_t\|_{L^2}+\|\nabla u_t\|_{L^2})\nonumber\\
&\leq \frac{1}{12}\|\nabla u_t\|_{L^2}^2+C\psi^\alpha\|\sqrt{\rho}u_t\|_{L^2}^2+C\psi^\alpha,\\
R_{12}&\leq C\int\rho|\nabla\Phi||u_t||\nabla u|dx\nonumber\\
&\leq C\|\rho\|_{L^\infty}^\frac{1}{2}\|\nabla\Phi\|_{L^\infty}\|\sqrt{\rho}u_t\|_{L^2}\|\nabla u\|_{L^2}\nonumber\\
&\leq C\psi^\alpha\|\sqrt{\rho}u_t\|_{L^2}^2+C\psi^\alpha.\label{c67}
\end{align}
Substituting \eqref{c57}-\eqref{c67} into \eqref{c54}, and we get
\begin{align}\label{p68}
&\frac{1}{2}\frac{d}{dt}\int\rho|u_t|^2dx
+\int\big((2\mu+\lambda)({\rm div}u_t)^2+\mu\omega_t^2+R^{-1}|u_t|^2\big)dx\nonumber\\
&\leq C\psi^\alpha(1+\|\sqrt{\rho}u_t\|_{L^2}^2+\|\nabla^2 u\|_{L^2}^2)
+C\psi^\alpha\|\bar{x}^\frac{a}{2}\nabla\eta\|_{L^2}^2+\frac{1}{2}\|\bar{x}^\frac{a}{2}\Delta\eta\|_{L^2}^2\nonumber\\
&\leq \frac{1}{2}\|\bar{x}^\frac{a}{2}\Delta\eta\|_{L^2}^2+
C\psi^\alpha\|\sqrt{\rho}u_t\|_{L^2}^2+C\psi^\alpha\|\bar{x}^\frac{a}{2}\nabla\eta\|_{L^2}^2
+C\|\eta_t\|_{L^2}^2+C\psi^\alpha,
\end{align}
where in the last inequality we have used \eqref{p48}.

Next, we should estimate $\|\bar{x}^\frac{a}{2}\Delta\eta\|_{L^2}^2$. Indeed, multiplying $\eqref{b1}_3$ by $\Delta\eta\bar{x}^a$,
integrating the result equation by parts yields that
\begin{align}\label{c69}
&\frac{1}{2}\Big(\int|\nabla\eta|^2\bar{x}^adx\Big)_t
+\int|\Delta\eta|^2\bar{x}^adx\nonumber\\
&=-\int\eta_t\cdot\nabla\eta\nabla\bar{x}^adx
+\int\nabla\cdot(\eta u-\eta\nabla\Phi)\cdot\Delta\eta\bar{x}^adx\nonumber\\
&=-\int(\Delta\eta-\nabla\cdot(\eta u-\eta\nabla\Phi))\cdot\nabla\eta\nabla\bar{x}^adx
+\int\nabla\cdot(\eta u-\eta\nabla\Phi)\cdot\Delta\eta\bar{x}^adx\nonumber\\
&\leq C\int|\nabla\eta||\Delta\eta||\nabla\bar{x}^a|dx+C\int|\nabla\eta|^2|u||\nabla\bar{x}^a|dx
+C\int\eta|\nabla\eta||\nabla u||\nabla\bar{x}^a|dx\nonumber\\
&\quad+C\int|\nabla\eta|^2|\nabla\Phi||\nabla\bar{x}^a|dx+
C\int\eta|\nabla\eta||\Delta\Phi||\nabla\bar{x}^a|dx
+C\int|\nabla u||\nabla\eta|^2\bar{x}^adx\nonumber\\
&\quad+C\int\eta|\nabla u||\Delta\eta|\bar{x}^adx+C\int|\nabla\eta||\nabla\Phi||\Delta\eta|\bar{x}^adx
+C\int\eta|\Delta\Phi||\Delta\eta|\bar{x}^adx\nonumber\\
&\triangleq \sum_{i=1}^9S_i.
\end{align}
Using the Gagliardo-Nirenberg inequality, \eqref{b15}, \eqref{p25}, \eqref{b4} and \eqref{b5}, we get
\begin{align}
S_1&\leq C\int\bar{x}^\frac{a}{2}|\nabla\eta|\bar{x}^\frac{a}{2}|\Delta\eta|\bar{x}^{-1}\log^{1+\sigma_0}(e+|x|^2)dx\nonumber\\
&\leq \varepsilon\|\bar{x}^\frac{a}{2}\Delta\eta\|_{L^2}^2+C\|\bar{x}^\frac{a}{2}\nabla\eta\|_{L^2}^2,\label{c70}\\
S_2&\leq C\int|\nabla\eta|^\frac{2a-1}{a}\bar{x}^\frac{2a-1}{2}|\nabla\eta|^\frac{1}{a}|u|\bar{x}^{-\frac{1}{4}}
\bar{x}^{-\frac{1}{4}}\log^{1+\sigma_0}(e+|x|^2)dx\nonumber\\
&\leq C\|\bar{x}^\frac{2a-1}{2}|\nabla\eta|^\frac{2a-1}{a}\|_{L^\frac{2a}{2a-1}}
\|u\bar{x}^{-\frac{1}{4}}\|_{L^{4a}}\||\nabla\eta|^\frac{1}{a}\|_{L^{4a}}\nonumber\\
&\leq C\psi^\alpha\|\bar{x}^\frac{a}{2}\nabla\eta\|_{L^2}^2+\varepsilon\|\bar{x}^\frac{a}{2}\Delta\eta\|_{L^2}^2,\\
S_3&\leq C\int\bar{x}^\frac{a}{2}\eta|\nabla u|\bar{x}^\frac{a}{2}|\nabla\eta|\bar{x}^{-1}\log^{1+\sigma_0}(e+|x|^2)dx\nonumber\\
&\leq C\|\bar{x}^\frac{a}{2}\eta\|_{L^4}\|\nabla u\|_{L^4}\|\bar{x}^\frac{a}{2}\nabla\eta\|_{L^2}\nonumber\\
&\leq C\|\bar{x}^\frac{a}{2}\eta\|_{L^4}^4+C\|\nabla u\|_{L^4}^4+C\|\bar{x}^\frac{a}{2}\nabla\eta\|_{L^2}^2\nonumber\\
&\leq C\|\bar{x}^\frac{a}{2}\eta\|_{L^2}^2(\|\bar{x}^\frac{a}{2}\nabla\eta\|_{L^2}^2
+\|\bar{x}^\frac{a}{2}\eta\|_{L^2}^2)+C\|\nabla u\|_{L^4}^4+C\|\bar{x}^\frac{a}{2}\nabla\eta\|_{L^2}^2\nonumber\\
&\leq C\psi^\alpha(\|\bar{x}^\frac{a}{2}\nabla\eta\|_{L^2}^2+\|\nabla^2 u\|_{L^2}^2),\\
S_4&\leq C\int|\nabla\eta|^2|\nabla\Phi|\bar{x}^a\bar{x}^{-1}\log^{1+\sigma_0}(e+|x|^2)dx\nonumber\\
&\leq C\|\nabla\Phi\|_{L^\infty}\|\bar{x}^\frac{a}{2}\nabla\eta\|_{L^2}^2,\\
S_5&\leq C\int \bar{x}^\frac{a}{2}\eta|\Delta\Phi|\bar{x}^\frac{a}{2}
|\nabla\eta|\bar{x}^{-1}\log^{1+\sigma_0}(e+|x|^2)dx\nonumber\\
&\leq C\|\Delta\Phi\|_{L^\infty}\|\bar{x}^\frac{a}{2}\eta\|_{L^2}^2
+C\|\bar{x}^\frac{a}{2}\nabla\eta\|_{L^2}^2\\
S_6&\leq C\|\nabla u\|_{L^2}\|\bar{x}^\frac{a}{2}\nabla\eta\|_{L^4}^2\nonumber\\
&\leq C\|\nabla u\|_{L^2}\|\bar{x}^\frac{a}{2}\nabla\eta\|_{L^2}(\|\bar{x}^\frac{a}{2}\nabla\eta\|_{L^2}
+\|\bar{x}^\frac{a}{2}\Delta\eta\|_{L^2})\nonumber\\
&\leq \varepsilon\|\bar{x}^\frac{a}{2}\Delta\eta\|_{L^2}
+C(\varepsilon)\psi^\alpha\|\bar{x}^\frac{a}{2}\nabla\eta\|_{L^2},\\
S_7&\leq C\|\bar{x}^\frac{a}{2}\Delta\eta\|_{L^2}
\|\bar{x}^\frac{a}{2}\eta\|_{L^4}\|\nabla u\|_{L^4}\nonumber\\
&\leq \varepsilon\|\bar{x}^\frac{a}{2}\Delta\eta\|_{L^2}^2
+C\|\bar{x}^\frac{a}{2}\eta\|_{L^2}\|\bar{x}^\frac{a}{2}\nabla\eta\|_{L^2}\|\nabla u\|_{L^2}\|\nabla u\|_{H^1}\nonumber\\
&\leq \varepsilon\|\bar{x}^\frac{a}{2}\Delta\eta\|_{L^2}^2
+C\psi^{-1}\|\nabla^2 u\|_{L^2}^2+C\psi^\alpha\|\bar{x}^\frac{a}{2}\nabla\eta\|_{L^2}^2,\\
S_8+S_9&\leq C(\|\nabla\Phi\|_{L^\infty}\|\bar{x}^\frac{a}{2}\nabla\eta\|_{L^2}
+\|\Delta\Phi\|_{L^\infty}\|\bar{x}^\frac{a}{2}\eta\|_{L^2})\|\bar{x}^\frac{a}{2}\Delta\eta\|_{L^2}\nonumber\\
&\leq \varepsilon\|\bar{x}^\frac{a}{2}\Delta\eta\|_{L^2}^2+C\|\bar{x}^\frac{a}{2}\nabla\eta\|_{L^2}^2
+C\psi^\alpha.\label{c77}
\end{align}
Substituting \eqref{c70}-\eqref{c77} into \eqref{c69} and choosing $\varepsilon$ suitably small lead to
\begin{align}\label{p78}
&\frac{1}{2}\frac{d}{dt}\|\bar{x}^\frac{a}{2}\nabla\eta\|_{L^2}^2+\|\bar{x}^\frac{a}{2}\Delta\eta\|_{L^2}^2\nonumber\\
& \quad \leq \varepsilon\|\bar{x}^\frac{a}{2}\Delta\eta\|_{L^2}^2+C\psi^\alpha\|\bar{x}^\frac{a}{2}\nabla\eta\|_{L^2}^2
+C\psi^{-1}\|\nabla^2 u\|_{L^2}^2+C\psi^\alpha.
\end{align}
Thus, multiplied  \eqref{p78} by $s$, together with Gronwall's inequality, we get
\begin{align}\label{p79}
\sup_{0\leq s\leq t}s\|\bar{x}^\frac{a}{2}\nabla\eta\|_{L^2}^2
+\int_0^ts\|\bar{x}^\frac{a}{2}\Delta\eta\|_{L^2}^2ds\leq C\exp\Big\{\int_0^t\psi^\alpha ds\Big\},
\end{align}
due to \eqref{p48} and \eqref{c12}.

Now, multiplying \eqref{p68} by $t$, we obtain \eqref{p52} using Gronwall's inequality and \eqref{p79}.
 The proof of Lemma \ref{lm3} is completed.
\begin{lem}\label{bx}~
Let $(\rho, u, \eta)$ and $T_1$ be as in Lemma \ref{xx}. Then, for all $t\in (0, T_1]$,
\begin{align}\label{d80}
\sup_{0\leq s\leq t}\|\bar{x}^a\rho\|_{L^1\cap H^1\cap W^{1, q}}
\leq \exp\Big\{C\exp\Big\{\int_0^t\psi^\alpha ds\Big\}\Big\}.
\end{align}
\end{lem}
{\bf Proof.}~ Notice that following the framework of Lemma 3.4 in \cite{LJ2} for proving an
estimate similar to \eqref{d80}, it suffices to verify the
following estimate:
\begin{equation}\label{c81}
\int_0^t\big(\|\nabla^2 u\|_{L^2\cap L^q}^{\frac{q+1}{q}}+s\|\nabla^2
u\|_{L^2\cap L^q}^2\big)ds\le  C\exp\Big\{C\int_0^t\psi^\alpha
ds\Big\}.
\end{equation}
In fact, on the one hand, it follows from \eqref{p48}, \eqref{p52} and  \eqref{b14} that
\begin{align}\label{c82}
&\int_0^t\big(\|\nabla^2 u\|_{L^2}^\frac{5}{3}+s\|\nabla^2 u\|_{L^2}^2\big)ds\nonumber\\
&\leq C\int_0^t(\psi^\alpha+\|\sqrt{\rho}u_t\|_{L^2}^2)ds+C\sup_{0\leq s\leq t}
\big(s\|\sqrt{\rho}u_t\|_{L^2}^2\big)\int_0^t\psi ds\nonumber\\
&\leq C\exp\Big\{C\int_0^t\psi^\alpha ds\Big\}.
\end{align}
On the other hand, choosing $p=q$ in \eqref{p47}, using the Gagliardo-Nirenberg inequality
gives
\begin{align}
\|\nabla^2 u\|_{L^q}&\leq C(\|\rho\dot{u}\|_{L^p}
+\|(\eta+\beta\rho)\nabla\Phi\|_{L^p}
+\|\nabla p_F\|_{L^p}+\|\nabla\eta\|_{L^p})\nonumber\\
&\leq C(\|\rho u_t\|_{L^q}+\|\rho u\|_{L^{2q}}\|\nabla u\|_{L^{2q}}+\|\eta\|_{L^p}
+\psi^\alpha+\|\nabla\eta\|_{L^p})\nonumber\\
&\leq C\|\rho u_t\|_{L^2}^\frac{2(q-1)}{q^2-2}\|\rho u_t\|_{L^2}^\frac{q^2-2q}{q^2-2}
+C\psi^\alpha(1+\|\nabla^2 u\|_{L^2}^{1-\frac{1}{q}}+\|\nabla^2\eta\|_{L^2}^{1-\frac{1}{q}})\nonumber\\
&\leq C\psi^\alpha(\|\sqrt{\rho}u_t\|_{L^2}^\frac{2(q-1)}{q^2-2}
\|\nabla u_t\|_{L^2}^\frac{q^2-2q}{q^2-2}+\|\sqrt{\rho}u_t\|_{L^2})
+C\psi^\alpha(1+\|\nabla^2 u\|_{L^2}^{1-\frac{1}{q}}+\|\nabla^2\eta\|_{L^2}^{1-\frac{1}{q}}).
\end{align}
This combined with \eqref{c82} \eqref{p52}, \eqref{b4}, and \eqref{b5}
\begin{align}
\int_0^t\|\nabla^2 u\|_{L^q}^\frac{q+1}{q}ds
&\leq C\int_0^t\psi^\alpha s^{-\frac{q+1}{2q}}(s\|\sqrt{\rho}u_t|_{L^2}^2)^\frac{q^2-1}{q(q^2-2)}
(s\|\nabla u_t\|_{L^2}^2)^\frac{(q-2)(q+1)}{2(q^2-2)}ds\nonumber\\
&\quad+C\int_0^t\psi^\alpha\|\sqrt{\rho}u_t\|_{L^2}^\frac{q+1}{q}ds
+C\int_0^t\Big(1+\|\nabla^2 u\|_{L^2}^\frac{q^2-1}{q^2}+\|\nabla^2\eta\|_{L^2}^\frac{q^2-1}{q^2}\Big)ds\nonumber\\
&\leq C\sup_{0\leq s\leq t}(s\|\sqrt{\rho}u_t\|_{L^2}^2)^\frac{q^2-1}{q(q^2-2)}
\int_0^t\psi^\alpha s^{-\frac{q+1}{2q}}(s\|\nabla u_t\|_{L^2}^2)^\frac{(q-2)(q+1)}{2(q^2-2)}ds\nonumber\\
&\quad+C\int_0^t\Big(\psi^\alpha+\|\sqrt{\rho} u_t\|_{L^2}^2+
\|\nabla^2 u\|_{L^2}^\frac{5}{3}+\|\nabla^2\eta\|_{L^2}^2\Big)ds\nonumber\\
&\leq C\exp\Big\{C\int_0^t\psi^\alpha ds\Big\}\Big[1+\int_0^t\Big(\psi^\alpha+s^{-\frac{q^3+q^2-2q-1}{q^3+q^2-2q}}
+s\|\nabla u_t\|_{L^2}^2\Big)ds\Big]\nonumber\\
&\leq C\exp\Big\{C\int_0^t\psi^\alpha ds\Big\},
\end{align}
and that
\begin{align}\label{c85}
\int_0^ts\|\nabla^2 u\|_{L^q}^2ds\leq C\exp\Big\{C\int_0^t\psi^\alpha ds\Big\}.
\end{align}
One thus obtains \eqref{c81} from \eqref{c82}-\label{c85} and completes the proof of lemma \ref{bx}.

Now, proposition \ref{mk} is a direct consequence of lemmas \ref{bx1}-\ref{bx}.
\\{\bf Proof of proposition \ref{mk}. }~It follows from
\eqref{d80}, \eqref{b4}, \eqref{b5}, and \eqref{d11}  and that
\begin{align*}
\psi(t)\leq \exp\Big\{C\exp\Big\{C\int_0^t\psi^\alpha ds\Big\}\Big\}.
\end{align*}
Standard arguments thus yield that for $M\triangleq e^{Ce}$ and $T_0\triangleq\min\{T_1, (CM^\alpha)^{-1}\}$,
\begin{align*}
\sup_{0\leq t\leq T_0}\psi(t)\leq M,
\end{align*}
which together with \eqref{p48}, \eqref{c81} and \eqref{c12}. The proof of Proposition \ref{mk}
is thus completed.
\begin{lem}\label{bx5}~
Let $(\rho, u, \eta)$ be a smooth solution to the initial-boundary-value problem \eqref{b1},
and $T_0$ is obtained in proposition \eqref{mk}, then we have
\begin{align}\label{d86}
\sup_{0\leq s\leq T_0}\big(s\|\eta_t\|_{L^2}^2+s\|\Delta\eta\|_{L^2}^2\big)+\int_0^{T_0}s\|\nabla\eta_t\|_{L^2}^2ds\leq C.
\end{align}
\end{lem}
{\bf Proof. }~
Differentiating $\eqref{b1}_3$ with respect to $t$ shows
\begin{align}\label{d87}
\eta_{tt}+\nabla\cdot(\eta_t u+\eta u_t-\eta_t\nabla\Phi)-\Delta\eta_t=0,
\end{align}
Multiplying \eqref{d87} by $\eta_t$ and then integrating equation over $B_R$, integrating by parts,
we have
\begin{align}\label{d88}
&\frac{1}{2}\frac{d}{dt}\int|\eta_t|^2dx+\int|\nabla\eta_t|^2dx\nonumber\\
&=\int(\eta_t u+\eta u_t-\eta_t\nabla\Phi)\cdot\nabla\eta_tdx\nonumber\\
&\leq C\int|\eta_t|^2|\nabla u|dx
+C\int\eta|u||\Delta\Phi||\nabla\eta_t|dx\nonumber\\
&\quad+C\int\eta|u_t||\nabla\eta_t|dx+C\int|\eta_t||\nabla\Phi||\nabla\eta_t|dx\triangleq\sum_{i=1}^4K_i.
\end{align}
Using the H\"older's inequality, Gagliardo-Nirenberg inequality, we have
\begin{align}
K_1&\leq C\|\eta_t\|_{L^4}^2\|\nabla u\|_{L^2}\nonumber\\
&\leq C\|\eta_t\|_{L^2}\|\nabla\eta_t\|_{L^2}\|\nabla u\|_{L^2}\nonumber\\
&\leq \varepsilon\|\nabla\eta_t\|_{L^2}^2+C(\varepsilon)\|\eta_t\|_{L^2},\label{d89}\\
K_2&\leq \int\eta\bar{x}^\frac{a}{2}|u|\bar{x}^\frac{a}{2}|\nabla\Phi||\nabla\eta_t|dx\nonumber\\
&\leq C\|\nabla\Phi\|_{L^\infty}\|\nabla\eta_t\|_{L^2}\|\bar{x}^\frac{a}{2}\eta\|_{L^4}
\|u\bar{x}^{-\frac{a}{2}}\|_{L^4}\nonumber\\
&\leq \varepsilon\|\nabla\eta_t\|_{L^2}^2
+C\|\bar{x}^\frac{a}{2}\nabla\eta\|_{L^2},\\
K_3&\leq C\|\nabla\eta_t\|_{L^2}\|\eta|u_t|\|_{L^2}\nonumber\\
&\leq \varepsilon\|\nabla\eta_t\|_{L^2}^2+C\int\eta\bar{x}^\frac{a}{2}\eta|u_t|^2\bar{x}^{-\frac{a}{2}}dx\nonumber\\
&\leq \varepsilon\|\nabla\eta_t\|_{L^2}^2+C\|\eta\|_{L^4}
\|u_t\bar{x}^{-\frac{a}{2}}\|_{L^8}^2\|\bar{x}^\frac{a}{2}\eta\|_{L^2}\nonumber\\
&\leq \varepsilon\|\nabla\eta_t\|_{L^2}^2+C(\varepsilon)\|\nabla u_t\|_{L^2}^2
+C(\varepsilon)\|\rho^\frac{1}{2}u_t\|_{L^2}^2,\\
K_4&\leq C\|\nabla\Phi\|_{L^\infty}\|\eta_t\|_{L^2}\|\nabla\eta_t\|_{L^2}\nonumber\\
&\leq  \varepsilon\|\nabla\eta_t\|_{L^2}^2+C(\varepsilon)\|\eta_t\|_{L^2}^2.\label{d92}
\end{align}
Now, putting \eqref{d89} and \eqref{d92} into \eqref{d88}, and multiplying the resulting inequality by $s$, we
have after choosing $\varepsilon$ suitably small that
\begin{align}
&\frac{d}{dt}\big(s\|\eta_t\|_{L^2}^2\big)+s\|\nabla\eta_t\|_{L^2}^2\nonumber\\
&\leq C(s\|\eta_t\|_{L^2}^2)+C(\|\eta_t\|_{L^2}^2+s\|\bar{x}^\frac{a}{2}\nabla\eta\|_{L^2}^2+s\|\nabla u_t\|_{L^2}
+s\|\rho^\frac{1}{2}u_t\|_{L^2}),
\end{align}
which together with Gronwall's inequality and \eqref{b14} yields that
\begin{align}
\sup_{0\leq s\leq T_0}\big(s\|\eta_t\|_{L^2}^2+s\|\Delta\eta\|_{L^2}^2\big)+\int_0^{T_0}s\|\nabla\eta_t\|_{L^2}^2ds\leq C.
\end{align}
The proof of lemma \ref{bx5} is completed.
\section{Proofs of theorems \ref{thm1}}

Let $(\rho_0, u_0, \eta_0)$ be as in Theorem \ref{thm1}.
For simplicity, assume that
\begin{align*}
\int_{\Bbb R^2}\rho_0 dx=1,
\end{align*}
which implies that there exists a positive constant $N_0$ such that
\begin{align}
\int_{B_{N_0}}\rho_0 dx\ge \frac{3}{4}\int_{\Bbb R^2}\rho_0dx=\frac{3}{4}.
\end{align}
We construct $\rho_0^R=\hat{\rho}_0^R+R^{-1}e^{-|x|^2}$ where
$0\leq\hat{\rho}_0^R\in C_0^\infty(\Bbb R^2)$ satisfies that
\begin{align}
\int_{B_{N_0}}\hat{\rho}_0^Rdx\ge\frac{1}{2},
\end{align}
and that
\begin{align}
\bar{x}^a\hat{\rho}_0^R\rightarrow\bar{x}^a\rho_0, \quad {\rm in}~L^1(\Bbb R^2)\cap H^1(\Bbb R^2)\cap W^{1, q}(\Bbb R^2)
\quad {\rm as}~R\rightarrow\infty.
\end{align}
Notice that $\bar{x}^\frac{a}{2}\eta_0\in L^2(\Bbb R^2)$ and $\bar{x}^\frac{a}{2}\nabla\eta_0\in L^2(\Bbb R^2)$,
choosing $\eta_0^R\in C_0^\infty(B_R)$ such that
\begin{align}\label{45}
\bar{x}^\frac{a}{2}\eta_0^R\rightarrow\bar{x}^\frac{a}{2}\eta_0, \quad\nabla\eta_0^R\rightarrow\nabla
\eta_0\quad {\rm in}~L^2(\Bbb R^2), \quad {\rm as}~ R\rightarrow\infty.
\end{align}
Since $\nabla u_0\in L^2(\Bbb R^2)$, choosing $v_i^R\in C_0^\infty(B_R)(i=1, 2)$ such that
\begin{align}
\lim_{R\rightarrow\infty}\|v_i^R-\partial_i u_0\|_{L^2(\Bbb R^2)}=0, \quad i=1,2,
\end{align}
and let smooth $u_0^R$ uniquely solve
\begin{equation}\label{4q}
\left\{
\begin{array}{ll}
\displaystyle -\Delta
u_0^R+R^{-1}u_0^R=-\rho_0^Ru_0^R+\sqrt{\rho_0^R}h^R-\partial_iv_i^R
&\mathrm{in}~B_R,\\
u_0^R=0 &\mathrm{on}~\partial B_R,
\end{array}
\right.
\end{equation}
where $h^R=(\sqrt{\rho_0}w_0^R)\ast j_{1/R}$ with the standard
mollifying kernel  $j_\delta$, $\delta>0$. Extend $u_0^R$ to
$\mathbb{R}^2$ by defining 0 outside $B_R$, and denote
$w_0^R\triangleq u_0^R\varphi_R$. By the same arguments as those for
the proof of Theorem 1.1 in \cite{LJ2}, we obtained that
\begin{equation}\label{6rr}
\lim_{R\rightarrow\infty}\Big(
\|\nabla(w_0^R-u_0)\|_{L^2(\mathbb{R}^2)}+ \|\sqrt{\rho_0^R}
w_0^R-\sqrt{\rho_0}u_0\|_{L^2(\mathbb{R}^2)}\Big)=0,
\end{equation}
where
\begin{align}\label{pp}
0\leq \varphi_R\leq 1, \quad \varphi_R(x)=1, ~{\rm if} ~|x|\leq R/2,
\quad |\nabla^k\varphi_R|\leq CR^{-k}(k=1, 2).
\end{align}
Then, in terms of lemma \ref{mn}, the initial-boundary value problem \eqref{b1}
 with the initial data $(\rho_0^R, u_0^R, \eta_0^R)$ has a classical solution
 $(\rho^R, u^R, \eta^R)$ on $B_R\times[0, T_R]$. Moreover, proposition \ref{mk}
 show that exists a $T_0$ independent of $R$ such that \eqref{3.3} and \eqref{d86} hold for
 $(\rho^R, u^R, \eta^R)$. By \eqref{3.3}, \eqref{d11}, \eqref{45}, \eqref{6rr}, and \eqref{d86}, after taking a subsequence, $(\rho^R, u^R, \eta^R)$
 locally and weakly (in the corresponding spaces) converges to a strong solution $(\rho, u, \eta)$ of \eqref{a1}-\eqref{ee} on $\Bbb R^2\times (0, T_0]$
satisfying \eqref{e7} and \eqref{e8}. The proof of the existence part of theorem \ref{thm1} is completed.

Next prove the uniqueness of the strong solutions. Take two strong solutions $(\rho_i, u_i, \eta_i)(i=1, 2)$
 sharing the same initial data with \eqref{e7} and \eqref{e8}, and let $\bar{\rho}=\rho_2-\rho_1, \bar{u}=u_2-u_1,
 \bar{\eta}=\eta_2-\eta_1$. Then,
\begin{equation}\label{t1}
\left\{
\begin{array}{ll}
\bar{\rho}_t+(u_2\cdot\nabla)\bar{\rho}+\bar{u}\cdot\nabla\rho_1+
\bar{\rho}\operatorname{div}u_2+\rho_1\operatorname{div}\bar{u}=0,\\
\rho_{1}\bar{u}_t+\rho_1u_1\cdot\nabla\bar{u}
+\nabla(p_F(\rho_2)-p_F(\rho_1))+\nabla(\eta_2-\eta_1)\\
\quad=\mu\Delta \bar{u}+(\mu+\lambda)\nabla\operatorname{div}\bar{u}
-\bar{\rho}(u_{2t}+u_2\cdot\nabla u_2)-\rho_1\bar{u}\cdot\nabla
u_2-(\bar{\eta}+\beta\bar{\rho})\nabla\Phi,\\
\bar{\eta}_t+\nabla\cdot(\bar{\eta}u_2-\bar{\eta}\nabla\Phi)+\nabla\cdot(\eta_2\bar{u})
-\Delta\bar{\eta}=0.
\end{array}
\right.
\end{equation}
for $(x,t)\in \mathbb{R}^2\times(0,T_0]$ with
\begin{equation} \label{00}\bar{\rho}(x, 0)=\bar{ u}(x, 0)=\bar{\eta}(x,
0)=0, \quad x\in \mathbb{R}^2.
\end{equation}

Firstly, multiply $\eqref{t1}_1$ by $2\bar{\rho}\bar{x}^{2r}$ and
integrate by parts. Similar to  the inequality (5.32) in
\cite {LJ2}, we get that
\begin{equation}\label{mh1}
\|\bar{\rho}\bar{x}^r\|_{L^2}\leq C\int_0^t(\|\nabla
\bar{u}\|_{L^2}+\|\sqrt{\rho_1}\bar{u}\|_{L^2})ds,\quad t\in
(0,T_0],
\end{equation}
where $r\in (1, \tilde{a})$ with $\tilde{a}=\min\{a, 2\}$.

Secondly, multiplying $\eqref{t1}_2$ by $\bar{u}$ and integrating
by parts lead to
\begin{align}\label{gg1}
&\frac{1}{2}\frac{d}{dt}\int\rho_1|\bar{u}|^2dx
+\int(2\mu+\lambda)|{\rm div}\bar{u}|^2
+\mu|\omega|^2dx\nonumber\\
&=-\int\bar{\rho}(u_{2t}+u_2\cdot\nabla
u_2)\cdot\bar{u}dx-\int\rho_1\bar{u}\cdot\nabla u_2\cdot\bar{u}dx\nonumber\\
&\quad+
\int(p_F(\rho_2)-p_F(\rho_1)){\rm div}\bar{u}dx+\int\bar{\eta} {\rm div}\bar{u}dx
+\int(\bar{\eta}+\beta\bar{\rho})\nabla\Phi\bar{u}dx\nonumber\\
&\leq C\|\nabla u_2\|_{L^\infty}\int\rho_1|\bar{u}|^2dx+
C\int|\bar{\rho}||\bar{u}|(|u_{2t}|+|u_2||\nabla
u_2|)dx\nonumber\\
&\quad+C\|p_F(\rho_2)-p_F(\rho_1)\|_{L^2}\|{\rm div}\bar{u}\|_{L^2}
+C\|\bar{\eta}\|_{L^2}\|{\rm div}\bar{u}\|_{L^2}\nonumber\\
&\quad+C\int\bar{\eta}\bar{x}^\frac{r}{2}|\bar{u}|\bar{x}^{-\frac{r}{2}}|\nabla\Phi| dx
+C\int\bar{\rho}\bar{x}^\frac{r}{2}|\bar{u}|\bar{x}^{-\frac{r}{2}}|\nabla\Phi| dx\nonumber\\
&\triangleq C\|\nabla u_2\|_{L^\infty}\int\rho_1|\bar{u}|^2dx+\sum_{i=1}^5Q_i.
\end{align}
Just like (\ref{mh1}), it has been obtained via (5.33) and (5.36) in
\cite{LJ2} that
\begin{align}\label{x29}
Q_1+Q_2&\leq C(\varepsilon)(1+t\|\nabla u_{2t}\|_{L^2}^2+t\|\nabla^2
u_2\|_{L^q}^2)\int_0^t(\|\nabla\bar{u}\|_{L^2}^2+
\|\sqrt{\rho_1}\bar{u}\|_{L^2}^2)ds\nonumber\\
&\quad+\varepsilon(\|\sqrt{\rho_1}\bar{u}\|_{L^2}^2+\|\nabla\bar{u}\|_{L^2}^2).
\end{align}
With the Cauchy inequality and \eqref{3.3}, \eqref{3.10}, and \eqref{d86}, we have
\begin{align}
\sum_{i=3}^5Q_i&\leq \varepsilon\|\nabla\bar{u}\|_{L^2}^2+C(\varepsilon)\|\bar{\eta}\|_{L^2}^2
+C\|\nabla\Phi\|_{L^4}\|\bar{x}^\frac{r}{2}\bar{\eta}\|_{L^2}\|\bar{u}\bar{x}^{-\frac{r}{2}}\|_{L^4}\nonumber\\
&\quad+C\|\nabla\Phi\|_{L^4}\|\bar{x}^\frac{r}{2}\bar{\rho}\|_{L^2}\|\bar{u}\bar{x}^{-\frac{r}{2}}\|_{L^4}
\nonumber\\
&\leq \varepsilon\|\nabla\bar{u}\|_{L^2}^2+C(\varepsilon)\|\bar{\eta}\|_{L^2}^2
+C\|\bar{x}^\frac{r}{2}\bar{\eta}\|_{L^2}(\|\sqrt{\rho_1}\bar{u}\|_{L^2}+\|\nabla\bar{u}\|_{L^2})\nonumber\\
&\quad+C(\|\sqrt{\rho}\bar{u}\|_{L^2}+\|\nabla\bar{u}\|_{L^2})\int_0^t
(\|\nabla\bar{u}\|_{L^2}+\|\sqrt{\rho_1}\bar{u}\|_{L^2})ds\nonumber\\
&\leq \varepsilon(\|\nabla\bar{u}\|_{L^2}^2+\|\sqrt{\rho_1}\bar{u}\|_{L^2}^2)+C(\varepsilon)\|\bar{\eta}\|_{L^2}^2
+C(\varepsilon)\|\bar{x}^\frac{r}{2}\bar{\eta}\|_{L^2}^2\nonumber\\
&\quad+\int_0^t
(\|\nabla\bar{u}\|_{L^2}^2+\|\sqrt{\rho_1}\bar{u}\|_{L^2}^2)ds.
\end{align}
It remains to estimate the $\|\bar{\eta}\|_{L^2}$ and $\|\bar{x}^\frac{r}{2}\bar{\eta}\|_{L^2}$. In fact,
multiplying $\eqref{a1}_3$ by $\bar{\eta}$ and integrating
\begin{align}
&\frac{1}{2}\frac{d}{dt}\int|\bar{\eta}|^2dx+\|\nabla\bar{\eta}\|_{L^2}^2\nonumber\\
&=-\int\bar{\eta}u_2\cdot\nabla\bar{\eta}dx-\int\bar{\eta}^2{\rm div}u_2dx
+\int\bar{\eta}\nabla\bar{\eta}\cdot\nabla\Phi dx\nonumber\\
&\quad+\int\bar{\eta}^2\Delta\Phi dx
-\int\bar{\eta}\bar{u}\cdot\nabla\eta_2dx-\int\eta_2{\rm div}\bar{u}\bar{\eta}dx\nonumber\\
&\triangleq \sum_{i=1}^6II_i.
\end{align}
Using the H\"older's inequality, Gagliardo-Nirenberg inequality, we have
\begin{align}
II_1&\leq C\int\bar{x}^\frac{r}{2}\bar{\eta}\bar{x}^{-\frac{r}{2}}|u_2||\nabla\bar{\eta}|dx\nonumber\\
&\leq C\|\nabla\bar{\eta}\|_{L^2}\|\bar{x}^\frac{r}{2}\bar{\eta}\|_{L^4}\|u_2\bar{x}^{-\frac{r}{2}}\|_{L^4}\nonumber\\
&\leq \varepsilon\|\nabla\bar{\eta}\|_{L^2}^2+
C(\varepsilon)\|\bar{x}^\frac{r}{2}\bar{\eta}\|_{L^2}^2+\delta\|\bar{x}^\frac{r}{2}\nabla\bar{\eta}\|_{L^2}^2,\\
\sum_{i=2}^4II_i&\leq C\|\nabla u_2\|_{L^2}\|\bar{\eta}\|_{L^4}^2+C\|\Delta\Phi\|_{L^\infty}\|\eta\|_{L^2}^2\nonumber\\
&\quad+\|\nabla\Phi\|_{L^\infty}\|\bar{\eta}\|_{L^2}
\|\nabla\bar{\eta}\|_{L^2}\nonumber\\
&\leq \varepsilon\|\nabla\bar{\eta}\|_{L^2}^2+C(\varepsilon)\|\bar{\eta}\|_{L^2}^2
+C,\\
II_5&\leq C\int\bar{x}^\frac{r}{2}\bar{\eta}\bar{x}^{-\frac{r}{2}}|\bar{u}||\nabla\eta_2|dx\nonumber\\
&\leq C\|\nabla\eta_2\|_{L^2}\|\bar{u}\bar{x}^{-\frac{r}{2}}\|_{L^4}\|\bar{x}^\frac{r}{2}\bar{\eta}\|_{L^4}\nonumber\\
&\leq \varepsilon(\|\sqrt{\rho_1}\bar{u}\|_{L^2}^2+\|\nabla\bar{u}\|_{L^2}^2)+
C(\varepsilon)\|\bar{x}^\frac{r}{2}\bar{\eta}\|_{L^2}^2+\varepsilon\|\bar{x}^\frac{r}{2}\nabla\bar{\eta}\|_{L^2}^2,\\
II_6&\leq C\|\nabla\bar{u}\|_{L^2}\|\eta_2\|_{L^4}\|\bar{\eta}\|_{L^4}\nonumber\\
&\leq \varepsilon\|\nabla\bar{u}\|_{L^2}^2+\delta\|\nabla\bar{\eta}\|_{L^2}^2+C(\varepsilon,\delta)\|\bar{\eta}\|_{L^2}^2.
\end{align}
Moreover, multiplying $\eqref{a1}_3$ by $\bar{x}^r\bar{\eta}$ and integrating by parts yield
\begin{align}\label{t36}
&\frac{1}{2}\frac{d}{dt}\int|\bar{\eta}|^2\bar{x}^rdx+\int|\nabla\bar{\eta}|^2\bar{x}^rdx\nonumber\\
&=\int\bar{\eta}^2\Delta\bar{x}^rdx-\int u_2\cdot\nabla\bar{\eta}\bar{\eta}\bar{x}^rdx
-\int\bar{\eta}{\rm div}u_2\bar{\eta}\bar{x}^rdx\nonumber\\
&\quad+\int\nabla\bar{\eta}\nabla\cdot\Phi\bar{\eta}\bar{x}^rdx+\int\bar{\eta}^2\bar{x}^r\Delta\Phi dx
+\int\bar{u}\eta_2\nabla\bar{\eta}\bar{x}^rdx\nonumber\\
&\quad+\int\bar{u}\eta_2\bar{\eta}\nabla\bar{x}^rdx\triangleq\sum_{i=1}^7III_i.
\end{align}
For the term $III_i(i=1,\cdots, 7)$ on the right hand side of \eqref{t36}, we get that
\begin{align}
III_1&\leq C\int|\bar{\eta}|^2\bar{x}^r\bar{x}^{-2}\log^{2(1+\sigma_0)}(e+|x|^2)dx\leq
C\int|\bar{\eta}|^2\bar{x}^rdx,\\
III_2+III_3&=-2\int\bar{\eta}^2{\rm div}u_2\bar{x}^rdx-\int\bar{\eta}^2u_2\nabla\bar{x}^rdx\nonumber\\
&\leq C\|\nabla u_2\|_{L^2}\|\bar{x}^\frac{a}{2}\bar{\eta}\|_{L^4}^2
+C\|\bar{x}^\frac{r}{2}\bar{\eta}\|_{L^4}
\|\bar{x}^\frac{r}{2}\bar{\eta}\|_{L^2}\|u_2\bar{x}^{-\frac{3}{4}}\|_{L^4}\nonumber\\
&\leq \varepsilon\|\bar{x}^\frac{r}{2}\nabla\bar{\eta}\|_{L^2}^2
+C(\varepsilon)\|\bar{x}^\frac{r}{2}\bar{\eta}\|_{L^2}
+C,\\
III_4+III_5&\leq \|\nabla\Phi\|_{L^\infty}\|\bar{x}^\frac{r}{2}
\bar{\eta}\|_{L^2}\|\bar{x}^\frac{r}{2}\nabla\bar{\eta}\|_{L^2}+C\|\Delta\Phi\|_{L^\infty}\|\bar{x}^\frac{r}{2}\bar{\eta}\|_{L^2}\nonumber\\
&\leq \varepsilon\|\bar{x}^\frac{r}{2}\nabla\bar{\eta}\|_{L^2}^2
+C(\varepsilon)\|\bar{x}^\frac{r}{2}\bar{\eta}\|_{L^2}
+C,\\
III_6&\leq \int\bar{u}\bar{x}^{-\frac{b}{2}}
\eta_2\bar{x}^\frac{b+r}{2}\nabla\bar{\eta}\bar{x}^\frac{r}{2}dx\nonumber\\
&\leq \|\bar{x}^\frac{r}{2}\nabla\bar{\eta}\|_{L^2}\|\bar{u}\bar{x}^{-\frac{b}{2}}\|_{L^4}
\|\bar{x}^\frac{a}{2}\eta_2\|_{L^4}\nonumber\\
&\leq \varepsilon\|\bar{x}^\frac{r}{2}\nabla\bar{\eta}\|_{L^2}^2
+\varepsilon(\|\sqrt{\rho_1}\bar{u}\|_{L^2}^2+\|\nabla\bar{u}\|_{L^2}^2)+C,\\
III_7&\leq C\int\bar{u}\bar{x}^{-\frac{3}{4}}\eta_2\bar{x}^\frac{r}{2}
\bar{\eta}\bar{x}^\frac{r}{2}\bar{x}^{-\frac{1}{4}}\log^{1+\sigma_0}(e+|x|^2)dx\nonumber\\
&\leq C\|\bar{x}^\frac{r}{2}\bar{\eta}\|_{L^2}\|\bar{u}\bar{x}^{-\frac{3}{4}}\|_{L^4}\|\bar{x}^\frac{r}{2}\eta_2\|_{L^4}\nonumber\\
&\leq \varepsilon(\|\sqrt{\rho_1}\bar{u}\|_{L^2}^2+\|\nabla\bar{u}\|_{L^2}^2)
+C(\varepsilon)\|\bar{x}^\frac{r}{2}\bar{\eta}\|_{L^2}^2,\label{x41}
\end{align}
where $b+r<\tilde{a}$.

Denoting
\begin{align}
G(t)\triangleq \|\sqrt{\rho_1}\bar{u}\|_{L^2}^2+\|\bar{\eta}\|_{L^2}+\|\bar{x}^\frac{r}{2}\bar{\eta}\|_{L^2}
+\int_0^t(\|\nabla\bar{u}\|_{L^2}^2+\|\sqrt{\rho_1}\bar{u}\|_{L^2}^2)ds,
\end{align}
with all these estimates \eqref{x29}-\eqref{x41}, choosing $\varepsilon, \delta $ suitably small lead to
\begin{align}
G'(t)\leq C(1+\|\nabla u_2\|_{L^\infty}+t\|\nabla u_{2t}\|_{L^2}^2+t\|\nabla^2
u_2\|_{L^q}^2)G(t),
\end{align}
which together with Gronwall's inequality, and \eqref{e7} yields $G(t)=0$. Hence, $\bar{u}(x, t)=0$ and $\bar{\eta}(x, t)=0$
for almost everywhere $(x, t)\in \Bbb R^2\times (0, T_0)$. Then, one can deduce from \eqref{mh1} that $\bar{\rho}=0$
for almost everywhere $(x, t)\in \Bbb R^2\times (0, T_0)$. The proof of theorem \ref{thm1} is completed.

{\small
\end{document}